\documentclass
[
    a4paper,
    DIV=11,
]
{scrartcl}
\usepackage{graphicx}

\usepackage{subfig}
\usepackage[labelfont={bf},format={plain}]{caption}

\usepackage[a4paper, left=3cm, right=3cm, top=2.5cm, bottom=2.5cm]{geometry}

\usepackage{standalone}									% for input of network files
\usepackage{tikz}										% for input of network files
\usetikzlibrary{arrows}
\usetikzlibrary{arrows.meta}
\usetikzlibrary{shapes.geometric}

\usepackage{authblk}
\usepackage{amsmath,amssymb,amsfonts,amsthm}
\usepackage[pdffitwindow=false,
            plainpages=false,
            pdfpagelabels=true,
            pdfpagemode=UseOutlines,
            pdfpagelayout=SinglePage,
            bookmarks=false,
            colorlinks=true,
            hyperfootnotes=false,
            linkcolor=blue,
            urlcolor=blue!30!black,
            citecolor=green!50!black]{hyperref}
\usepackage[nameinlink]{cleveref}

\newtheorem{theorem}{Theorem}[section]
\newtheorem{corollary}[theorem]{Corollary}

\newtheorem{proposition}[theorem]{Proposition}
\newtheorem{definition}[theorem]{Definition}
\theoremstyle{definition}
\newtheorem{example}[theorem]{Example}
\newtheorem{remark}[theorem]{Remark}

\crefname{subsection}{\textup{subsection}}{\textup{subsections}}
\Crefname{subsection}{\textup{Subsection}}{\textup{Subsections}}
\crefname{proposition}{\textup{proposition}}{\textup{propositions}}
\Crefname{proposition}{\textup{Proposition}}{\textup{Propositions}}

\usepackage{enumitem}

\newcommand{\N}{\mathbb{N}}                                     % natural numbers
\newcommand{\R}{\mathbb{R}}  									% real numbers
\newcommand{\C}{\mathbb{C}}                                     % complex numbers

\newcommand{\NN}{\mathcal{N}}
\newcommand{\T}[1]{\widetilde{#1}}

\newcommand\set[1]{\left\{\,#1\,\right\}}
\newcommand\with{\ \vrule\ }
\providecommand{\norm}[1]{\left\lVert #1 \right\rVert}          % norm
\providecommand{\abs}[1]{\left\lvert #1 \right\rvert}
\DeclareMathOperator{\Span}{span}

\DeclareMathOperator{\normi}{norm_i}

\newcommand{\Syn}{V_0}
\DeclareMathOperator{\circulant}{circ}

\newcommand{\gtm}{\textsc{Go-To-The-Middle}}
\newcommand{\nbug}{\textsc{$N$-bug}}
\newcommand{\gta}{\textsc{Go-To-The-Average}}
\newcommand{\gtg}{\textsc{Go-To-The-Center-Of-Gravity}}
\newcommand{\fsync}{\textsc{$\mathcal{F}$sync}}

\def\cld/{\textsc{CLD}}
\def\nscld/{\textsc{NSCLD}}
\def\cmas/{\textsc{CMAS}}
\def\oblot/{$\mathcal{OBLOT}$}

\def\LCM/{\textsc{LCM}}
\def\Look/{\textsc{Look}}
\def\Compute/{\textsc{Compute}}
\def\Move/{\textsc{Move}}

\newcommand{\n}{\text{\fontfamily{pzc}\selectfont n}}

\title{On the Dynamical Hierarchy in Gathering Protocols with Circulant Topologies\thanks{This work is funded by the Deutsche Forschungsgemeinschaft (DFG, German Research Foun\-dation) -- Project number 453112019.}}
\author{Raphael Gerlach, Sören von der Gracht and Michael Dellnitz
\thanks{The authors are with the Institute of Mathematics, Paderborn University, Warburger Str. 100, 33098 Paderborn, Germany (e-mail: raphael.gerlach@uni-paderborn.de, soeren.von.der.gracht@uni-paderborn.de, dellnitz@uni-paderborn.de).}}
\date{}

\begin{document}

\maketitle

\begin{abstract}
In this article, we investigate the convergence behavior of two classes of gathering protocols with fixed circulant topologies using tools from dynamical systems.
Given a fixed number of mobile entities moving in the Euclidean plane, we model a gathering protocol as a system of (linear) ordinary differential equations whose equilibria are exactly all possible gathering points.
Then, for a circulant topology we derive a decomposition of the state space into stable invariant subspaces with different convergence rates by utilizing tools from dynamical systems theory. It turns out, that this decomposition is identical for every linear circulant gathering protocol, whereas only the convergence rates depend on the weights in interaction graph itself. 
In the second part, we consider a normalized nonlinear version of the equation of motion that is obtained by scaling the speed of each entity. Again, we find a similar decomposition of the state space that is based on our findings in the linear case. 
Finally, we also consider visibility preservation properties of the two classes of system.

\medskip
\end{abstract}


\vspace{.3cm}

\section{Introduction}
This article applies dynamical systems theory to the analysis of swarms of mobile robots, which are widely studied under the headline of \emph{distributed computing} (e.g. \cite{DBLP:series/lncs/FlocchiniPS19,H18}).
This precise and robust mathematical framework can be exploited to formally prove statements about the collective behavior but also to inform design choices for algorithms,
which recently has been exemplified in the investigation of symmetry properties of robot swarms performing \textsc{Near-Gathering} \cite{GGHHK25}.

\textbf{Gathering:}
The problem of interest in this article is the \emph{gathering problem} in which the robots are supposed to converge to a single, not predefined, point. 
The only capabilities the robots have are observing other robots' positions, performing computations in their local memory, and moving.
The strategy that each robot pursues is called a \emph{protocol}.
While most studies of the gathering problem employ discrete time 
(e.g. \cite{DBLP:journals/tcs/CastenowFHJH20,Castenow.2023,DBLP:journals/siamcomp/CohenP05,DBLP:conf/spaa/DegenerKLHPW11,DBLP:series/lncs/Flocchini19}), some protocols using continuous time have been proposed as well (\cite{DBLP:journals/topc/DegenerKKH15,DBLP:series/lncs/KlingH19}). The latter framework is the focus of this manuscript. 
Gathering is an instance of famous problems known under terms such as \emph{synchronization, agreement problems, rendezvous problems, consensus algorithms}, and likely many more.
For example, in synchronization literature one commonly studies Laplacian (or Laplacian like) dynamics of coupled oscillators. Very general internal dynamics and coupling functions are possible and synchronization is not restricted to equilibria (see e.g.~\cite{Dorfler.2014,Heagy.1994,Pikovskij.2003} and references therein).
Very powerful synchronization results have also been derived in the control theory literature under the headline of \emph{distributed consensus control/algorithms} 
or \emph{asymptotic agreement problems}. Asymptotic consensus was investigated  for abstract variables, time varying interaction topologies, communication delay, and even defective agents (e.g. \cite{Chen.2013,Moreau.2004,OlfatiSaber.2003b,OlfatiSaber.2003,OlfatiSaber.2004}).

\textbf{The Circulant Multi-Agent System Model (\cmas/):}
We consider a swarm of $N\in \N$ autonomous, mobile robots modeled as deterministic points in the $2$-dimensional plane.
The robots are \emph{oblivious} (have no memory), disoriented (own local self-centered coordinate system without agreement on direction and chirality) and to have a fixed (ordered) list of \emph{interaction partners}. 
This list is realized as a fixed \emph{interaction graph}.
We assume the interaction graph, respectively its adjacency matrix, to be \emph{circulant} which allows all robots to execute the same algorithm (\emph{homogeneous}). This interaction topology can be seen as a generalization of a fixed closed chain.
We emphasize, that the robots are \emph{not} \emph{anonymous} as they can distinguish their interaction partners.
Interaction consists of the observation of the partner's positions which are used for the internal computations. 
These are not necessarily symmetric.
Finally, we assume that the robots have a limited, constant \emph{viewing range}, such that they cannot perceive anything beyond that range.
Robots move according to a \emph{time continuous} \emph{fully-synchronous} \fsync\ model, i.e, each robot continuously adjusts its speed and direction by synchronously performing \Look/-\Compute/-\Move/ (\LCM/) rounds permanently and instantaneously.

\textbf{Our Contribution:}
We propose two \emph{classes} of protocols called \emph{Circulant Laplacian Dynamics} and \emph{Nonlinear Scaling of Circulant Laplacian Dynamics}.
Protocols from the first class solve the gathering problem using the \cmas/ model
via convergence in infinite time. The scaled version facilitates finite time gathering as well.
For their respective analyses, we apply the mathematical theory of dynamical systems which allows for a profound understanding of the collective dynamics in terms of fine grained gathering rates.
We unveil a foliation of the space of all robots' positions into dynamically invariant subspaces in which the configurations gather with different speeds (\Cref{thm:decomp,thm:decomp_nonlinear}). This foliation is independent of the precise protocol. It allows to identify gathering rates for specific initial configurations and to decompose arbitrary initial configurations into components with different gathering rates (\Cref{cor:decomp}).
Moreover, we show that visibility is preserved under the dynamics (cf. \Cref{thr:visibility-circulant} and \Cref{prop:visibility-circulant-nonlinear}).

\begin{example}
	\label{ex:intro}
	Consider $N=7$ mobile robots $\set{0,\dotsc,6}$ running the \nbug{} protocol, i.e., if $z_i(t)\in\R^2$ is the position of robot $i$ at time $t$ the $i$-th robot is influenced by the $i+1$-th robot and will move towards it.
	\Cref{thm:decomp,thm:decomp_nonlinear} show that the configuration in \Cref{fig:intro}~(a) has the slowest gathering rate, while the configuration in \Cref{fig:intro}~(b) has the fastest gathering rate. Further intermediate configurations are unveiled as well. 
	\begin{figure}[!htb]
		\centering
		\subfloat[][\centering Slowest gathering configuration]{
			\begin{minipage}[t]{.45\linewidth}
				\centering
				\includegraphics[width=0.5\textwidth]{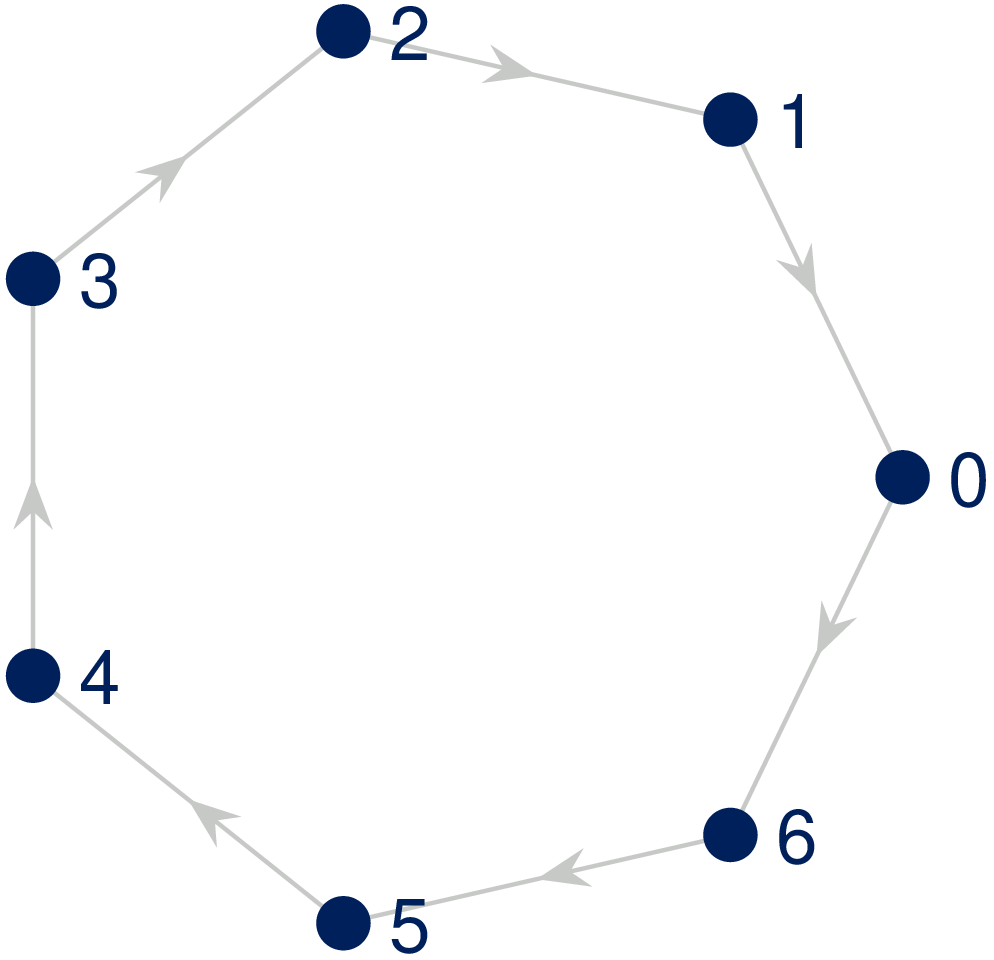}
			\end{minipage}
		}
		\hfil
		\subfloat[][\centering Fastest gathering configuration]{
			\begin{minipage}[t]{.45\linewidth}
				\centering
				\includegraphics[width=0.5\textwidth]{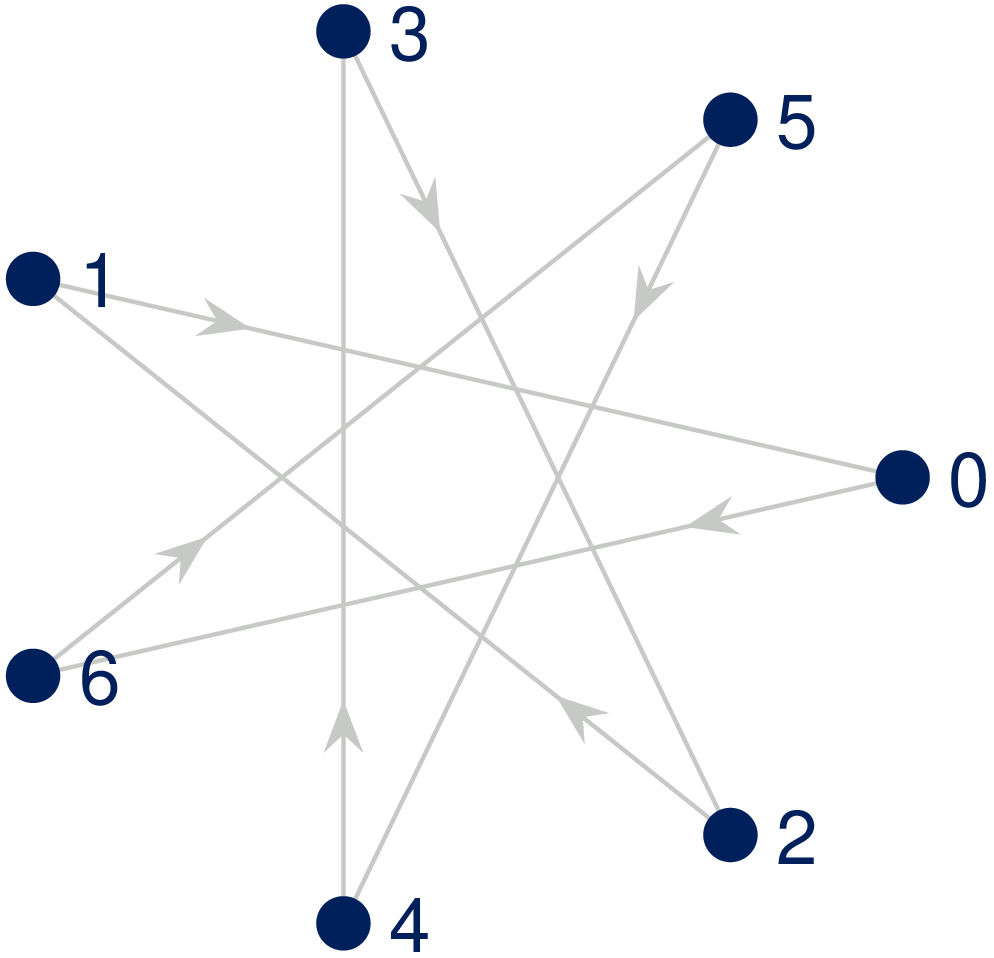}
			\end{minipage}
		}  
		\caption{Illustration of the gathering hierarchy exhibited by $N=7$ robots running the \nbug{} protocol. Blue dots indicate positions of each robot. Gray arrows symbolize the interaction relation, and \emph{not} the movement direction.}
		\label{fig:intro}
	\end{figure}
\end{example}

\textbf{Further Related Work:}
The fixed circulant interaction topology serves mainly as an academic example. Nonetheless, it can be used in the theoretical worst case run-time analysis of distributed strategies to form communication \emph{chains} between stationary relays (cf. \cite{Brandes.2013,Degener.2011b,Kling.2010,Kling.2011}).

Furthermore, we would like to mention that comparable results are known in the synchronization literature for coupled oscillators. There, dynamical hierarchy corresponds to \emph{normal modes} or \emph{eigenmodes} which can unveil communities that synchronize faster than others (see for example \cite{Arenas.2006,McGraw.2008}). Moreover, overall convergence speed to synchrony is commonly related to spectral properties of the Laplacian or system matrix, which has been done  in particular for circulant interaction topologies (e.g. \cite{Iqbal.2018,Irofti.2016,OlfatiSaber.2004,Rubido.2016}).
However, to the best of our knowledge a similarly detailed exploitation of the full spectral decomposition of phase space and the interpretation in terms of a dynamical hierarchy of configurations in positional variables has not been done before. 

\textbf{Outline:}
The article is organized as follows: In \Cref{sec:preliminaries}, we gather notation to formalize the underlying interaction structure by a weighted graph.
In \Cref{subsec:dynamics}, we propose two classes of protocols (\emph{Circulant Laplacian Dynamics (\cld/)} and \emph{Nonlinear Scaling of Circulant Laplacian Dynamics (\nscld/})).
We also present the mathematical dynamical systems framework for their analysis. The main dynamical analysis is carried out in \Cref{sec:analysis_lin} for the linear class and in \Cref{sec:analysis_nonlin} for the nonlinear class. In \Cref{subsec:hierachy_lin}, the hierarchy of different gathering rates is proven, which is extended in \Cref{sec:analysis_nonlin} to the nonlinear setting. Visibility preservation is shown in \Cref{subsec:preservation_lin} and \Cref{sec:analysis_nonlin}. Finally, in \Cref{sec:conclusion} we end this article with a conclusion and a brief outlook on future research. Details of technical mathematical proofs can be found in \Cref{app:apendix}.

\section{Circulant Interaction Topology}
\label{sec:preliminaries}
\label{subsec:interaction}
We begin by introducing some notation from algebraic graph theory.
As we assume the interaction structure to be fixed for all times independent of the robots' positions, we can encode interactions between robots by a (directed) \emph{graph} $G=(V,E)$ that we call \emph{interaction graph}. The \emph{vertices} are the robots $V = \set{0, \dotsc, N-1}$ and an \emph{edge} $e=(j,i)\in E \subseteq V\times V$ represents the fact that robot $i$ uses the position of robot $j$ in its internal computations.
We assume the interaction graph to be \emph{weighted}: each edge $e=(j,i)\in E$ is assigned the \emph{weight} $w_{i,j} \in \R$ and $w_{i,j}=0$ if and only if  $e=(j,i)\notin E$. 
A negative weight corresponds to the partner's position \emph{reflected at the origin} and would force a robot to intentionally \emph{move away} from that partner (cf. \Cref{sec:protocols} below). Hence, we consider non-negative weights only from now on.
A commonly considered special case is when the interaction graph is \emph{symmetric} or \emph{undirected}, i.e, the weight matrix $W\in \R^{N\times N}$ is symmetric.
The interaction graph is said to be \emph{connected} (or \emph{weakly connected}) if for any two vertices ${j,i \in\set{0,\dotsc,N-1}}$ there is an \emph{undirected path} from $j$ to $i$. If there is always a \emph{directed path} the graph is said to be \emph{strongly connected}. 

The interaction graph is \emph{circulant} if and only if there are integers $0\le s_1 < \dotsb < s_k \le N-1$ such that
$(j,i) \in E$ if and only if $j = i+s_r \bmod N$  for an $r\in\set{1,\dotsc,k}$.
Therein $i+s \bmod N = i+s$ if $i+s< N$ and $i+s-N$ if  $i+s \ge N$.
Circulant interaction structure implies that the weight matrix $W\in \R^{N\times N}$ is a \emph{circulant matrix} of the form
\begin{equation*}\label{eq:circulant}
	W=\circulant(w_0,w_1,\dotsc,w_{N-1}) := \left(\begin{smallmatrix}w_{0}&w_{1}&\cdots &w_{N-2}&w_{N-1}\\
		w_{N-1}&w_{0}&w_{1}&&w_{N-2}\\
		\vdots &w_{N-1}&w_{0}&\ddots &\vdots \\
		w_{2}&&\ddots &\ddots &w_{1}\\
		w_{1}&w_{2}&\cdots &w_{N-1}&w_{0}
	\end{smallmatrix}\right) \in \R^{N\times N}.
\end{equation*}
We say the vector $w=(w_0,\dotsc,w_{N-1})\in \R^N$ \emph{generates} the circulant matrix $W$. 
For a symmetric circulant weight matrix $W$, this implies
$w_{N-i} = w_i i=1,\dotsc,N-1)$
and $W$ has fewer degrees of freedom, i.e, $W$ is determined by only $\lfloor \frac{N}{2}\rfloor +1$ elements $w_0,\dotsc,w_{\lfloor\frac{N}{2}\rfloor}$.
For a circulant graph, weak and strong connectivity coincide (cf.~\cite[Corollary 1]{vanDoorn.1986}). Hence, we drop the distinction.

\section{Two classes of protocols}
\label{subsec:dynamics}
\label{sec:protocols}
We consider two classes of protocols: \emph{Circulant Laplacian Dynamics (\cld/)} and \emph{Nonlinear Scaling of Circulant Laplacian Dynamics (\nscld/)}. In both cases, a robot~$i$ continuously and instantaneously adapts its velocity vector according to the
following
rules corresponding to the \LCM/ phases: Robot $i$
\begin{description}[font=\normalfont\scshape]
	\item[Look.] observes the positions of the robots $j$ such that $e=(j,i) \in E$ in relative coordinates within its vision range,
	\item[Compute.] computes a target point as a linear combination of the observed positions where the coefficients are given by the weights of the interaction graph, and
	\item[Move.] adapts its velocity vector as the connecting vector between its own position and the target point. In an \nscld/ protocol, this vector is rescaled.
\end{description}

For the dynamical analysis, we take the stance of an external observer describing the dynamics in global coordinates. 
We denote by $z_i(t) = (x_i(t),y_i(t)) \in \R^2$ the position of robot $i \in \{0,\dotsc,N-1\}$ at time $t \in \R$. Typically, the time argument will be omitted for notational brevity. 
The collection of all robots' positions is called a \emph{configuration}.
The \LCM/ steps of a \cld/ protocol described above are mathematically modeled by the \emph{linear} ordinary differential equations
\begin{equation}
	\label{eq:f_lin_circ}
	\dot{z}_i = -z_i + \sum_{j=0}^{N-1} w_j z_{i+j},
\end{equation}
where robot indices are counted $\bmod~N$, i.e., $z_{(N-1)+1}=z_0$ and $z_{0-1}=z_{N-1}$.

\cld/ protocols satisfying the linear equations \eqref{eq:f_lin_circ} are idealized. They do not generally prevent unbounded velocity and realize gathering only in the limit $t\to\infty$, as we will see below in \Cref{sec:analysis_lin}. Both restrictions are at least problematic for the robot model.
However, in commonly investigated protocols robots compute their \emph{direction} as a linear combination of the (relative) positions of their neighbors (as in \eqref{eq:f_lin_circ}) and move with bounded speed until they are gathered. This requires rescaling the velocity vector in the \textsc{Move} step of a \nscld/ protocol. The modification is realized by a \emph{nonlinear} function $\NN:\R^2\to\R^2$ with $\NN(0) = 0$ of the form
\begin{align}\label{eq:N}
	\NN(x,y) = \n(\norm{(x,y)}_2)\begin{pmatrix}
		x\\y
	\end{pmatrix}
\end{align}
for some non-negative function $\n:\R \to \R_{\geq 0}$. This function \emph{scales} the length of the original velocity vector, and therefore the robot's speed, but does not change the direction.
Here, $\|\cdot\|$ denotes the Euclidean norm. Other norms are possible. In fact, the theoretical results below can be adapted without much effort. 
The mathematical model is adapted to the \emph{nonlinear} ordinary differential equations
\begin{equation}
	\label{eq:f_lin_normalized}
	\dot{z}_i  = \NN\left(-z_i + \sum_{j=0}^{N-1} w_j z_{i+j}\right).
\end{equation}
\noindent
A common approach is robots moving with constant speed one which requires
\begin{align}\label{eq:NN}
	\n(\norm{(x,y)}_2) = \norm{(x,y)}_2^{-1}.
\end{align}

\begin{remark}
	\label{rem:finite_viewing_range}
	Both mathematical formulations \eqref{eq:f_lin_circ} and \eqref{eq:f_lin_normalized} do not incorporate the finite viewing range of robots explicitly. Instead, we analyze the frameworks under the hypothetical assumption of unlimited viewing range and prove in hindsight that interaction partners do not lose sight of each other when they started suitably close. This approach is similar to commonly studied multi-agent systems in which the communication graph is all-to-all connected and communication failure corresponds to vanishing weights. Examples include flocking models \cite{Cucker.2007}, bounded confidence opinion dynamics \cite{Hegselmann.2002}, or consensus algorithms \cite{Moreau.2004}.
\end{remark}

In \Cref{sec:analysis_lin,sec:analysis_nonlin}, we derive conditions under which a protocol from the two classes \cld/ and \nscld/ solves the gathering problem.
To make the notion precise, we state.
\begin{definition}\label{def:gathering}
	A protocol is called \emph{gathering}, if the solution with any initial configuration converges to a \emph{gathering point} $z_i=z^*=(x^*,y^*) \in \R^2$ for all $i$ 
	for $t\to \infty$ and if the system is in equilibrium in any such gathering point, i.e., $\dot{z}_i = 0$ for all $i$.
	The subspace of all gathering configurations is denoted by 
	\[
	\Syn=\set{Z^*=(x^*, \dotsc,x^*, y^*,\dotsc,y^*) \with x^*,y^*\in\R} \subseteq \R^{2N}.
	\]
\end{definition}

\begin{example}\label{ex:running_ex}
	As a running example, we consider the \nbug{} problem \cite{Watton.1969}. Here, each robot~$i$ is only influenced by its first neighbor to the right, which yields a connected circulant interaction structure of the form $W=\circulant(0,w_1,0\dotsc,0)$ with ${w_1\ne0}$.
\end{example}

\section{Analysis of the Circulant Laplacian Dynamical Model}
\label{sec:analysis_lin}
\subsection{Gathering}
\label{subsec:cld-gathering}
In this section we study the convergence properties of \cld/ protocols using tools from dynamical systems theory. By \Cref{def:gathering} it can readily be seen, that a necessary condition for a linear strategy \eqref{eq:f_lin_circ} to be gathering is $z^*=\sum_{j=0}^{N-1} w_{j}z^*$, which is satisfied if and only if the weights satisfy
\begin{align}\label{eq:consistent}
	\sum_{j=0}^{N-1} w_{j} = 1.
\end{align}
In this situation, we say that the weights, the weight matrix, and the linear strategy are \emph{consistent}. Under this condition, the right hand side of \eqref{eq:f_lin_circ} is a linear combination of the \emph{relative} positions $z_i-z_j$. In particular, the \cmas/ model is able to perform computations of the form \eqref{eq:f_lin_circ}.

For the analysis we collect all robots' positions in a single vector $Z=(X,Y) \in \R^{2N}$, where $X = (x_0,\dotsc,x_{N-1}) \in \R^N$ denotes the $x$-coordinates of the $N$ robots and ${Y = (y_0,\dotsc,y_{N-1}) \in \R^N}$ the $y$-coordinates. 
This allows us to rewrite \eqref{eq:f_lin_circ} as
\begin{align}\label{eq:Z_lin}
	\dot{Z} = \left(-\mathbf{I}_{2N} + \mathbf{W}\right)Z 
	\text{ with } 	\mathbf{W} := \mathbf{I}_2 \otimes W =
	\begin{pmatrix}
		W & \vline & \mathbf{0}\\\hline
		\mathbf{0} & \vline & W
	\end{pmatrix} \in \R^{2N \times 2N},
\end{align}
where $\mathbf{I}_k \in \R^{k\times k}$ denotes the identity matrix and $\mathbf{I}_2\otimes W$ the Kronecker product.

The behavior of linear systems of ordinary differential equations is well understood. 
In the particular case of a circulant interaction graph, the classification of gathering protocols follows from Theorem~3.2 in \cite{Beard.2003}, Theorem~4 in \cite{OlfatiSaber.2004}, and consistency condition \eqref{eq:consistent} which we saw to be necessary for a gathering protocol. In fact, under the assumption of consistency, the system matrix $-\mathbf{I}_{2N} + \mathbf{W}$ in \eqref{eq:Z_lin} equals $\mathbf{I}_2 \otimes L$, where $L\in \R^{N\times N}$ is the graph Laplacian matrix of the weighted interaction graph so that these results apply.
\begin{theorem}
	\label{thr:gathering-circulant}
	For a \cld/ protocol with weight matrix $W=\circulant(w_0,\dotsc,w_{N-1})$, where $w_i\ge 0~(i=0,\dotsc,N-1)$, the following are equivalent:
	\begin{enumerate}[label = (\roman*)]
		\item The protocol is gathering.
		\item The interaction graph is connected and the weight matrix is consistent.
	\end{enumerate}
\end{theorem}
Seminal works as \cite{Beard.2003,Moreau.2004,OlfatiSaber.2003b,OlfatiSaber.2003,OlfatiSaber.2004} and many others have classified gathering protocols (consensus control strategies in their context) also in more general settings.

\begin{example}\label{ex:circ_gathering}
	Consider the \nbug{} problem introduced in the running \Cref{ex:running_ex}. For this protocol, \Cref{thr:gathering-circulant} yields that it is gathering if and only if its generating weight vector $w\in \R^N$ has the form ${w=(0,1,0,\dotsc,0)^T \in \R^N}$.
\end{example}

\subsection{Preservation of Visibility}
\label{subsec:preservation_lin}

So far in our analysis, we have ignored the fact that robots can only use local information specified by their finite viewing range $\mathcal{C}>0$.
While synchronization results under similar constraints are available (e.g. \cite{OlfatiSaber.2003b}), it remains largely unclear if visibility of different agents with a limited vision radius can be preserved under the dynamics.
However, we may prove that this issue does not occur with \cld/ protocols \emph{which solve the gathering problem}.
In fact, \Cref{thr:visibility-circulant} below shows that no robot loses sight of any of its interaction partners if initially all their interaction partners are within viewing range.
The proof crucially depends on the time continuity of the model. It uses the consistency of the weight matrix of a gathering circulant protocol to show that the two maximally distant robots cannot move further away from each other. Note that robots can leave the convex hull of the configuration, however, the swarms diameter shrinks sufficiently fast, such that loss of sight is impossible. For more details, we refer to \Cref{app:vis}.

\begin{theorem}
	\label{thr:visibility-circulant}
	Consider a gathering \cld/ protocol with circulant weight matrix $W\in \R^{N\times N}$.
	Let ${z_0(0),\dotsc,z_{N-1}(0)\in\R^2}$ be a valid initial configuration, that is,
	$\|z_i(0)-z_j(0)\| \le \mathcal{C}$ for all    $(j,i)\in E$.
	Then, visibility is preserved under the dynamics, i.e.,
	\[ \|z_i(t)-z_j(t)\| \le \mathcal{C} \quad\text{for all } t\geq 0 \text{ and } (j,i)\in E. \]
\end{theorem}

\subsection{Hierarchy of Convergence Rates}
\label{subsec:hierachy_lin}
\Cref{thr:gathering-circulant} provides a simple condition for \cld/ protocols to be gathering.
For general linear systems such as \eqref{eq:f_lin_circ} it is well known that spectral properties of the system matrix determine the longtime behavior.
For a circulant matrix all eigenvalues and eigenvectors can be analytically computed as follows (e.g. \cite{G05}).
\begin{proposition}\label{prop:circ_spectral}
	Let $W = \circulant(w_0,\dotsc,w_{N-1}) \in \R^{N\times N}$ be a circulant matrix. Then its eigenvectors $v_j\in \C$ ($j=0,\dotsc,N-1)$ are given by
	\begin{align}\label{eq:EV}
		v_j = (1, \omega^j, \omega^{2j},\dotsc,\omega^{(N-1)j})^T \in \C^N,
	\end{align}
	where $\omega = \exp{\left(\frac{2\pi \mathbf{i}}{N}\right)}$ is a primitive $N$-th root of unity and its eigenvalues are
	\begin{align}\label{eq:EW}
		\lambda_j = \sum_{i=0}^{N-1} w_i\omega^{ij} \in \C \text{ for } j=0,\dotsc,N-1.
	\end{align}
\end{proposition}

\begin{remark}\label{rem:circ_EW}\quad
	\begin{enumerate}[label = (\alph*)]
		\item For $j=0$, we obtain the eigenvector $v_0=(1,\dotsc,1)^T\in \R^N$ with corresponding (simple) eigenvalue $\lambda_0 = \sum_{i=0}^{N-1} w_i$, which has to be one for the protocol to be gathering according to \Cref{thr:gathering-circulant}.
		\item Since $W\in \R^{N\times N}$ has real-valued entries, we immediately obtain $\lambda_j = \overline{\lambda_{N-j}} \text{ for } j=1,\dotsc, N-1$.
		For even $N=2k$ this implies $\lambda_k\in \R$. Moreover, if $W \in \R^{N\times N}$ is symmetric, this reduces to $\lambda_j = \lambda_{N-j} \text{ for } j=1,\dotsc, N-1,$ and the eigenvalues $\lambda_1,\dotsc, \lambda_{\lfloor\frac{N-1}{2}\rfloor} \in \R$ are duplicated.
		\item By \eqref{eq:EV}, all eigenvectors $v_j\in \C^N$ are independent of the generating vector $w\in \R^N$, i.e., the given \cld/ protocol.
		In \Cref{fig:EV_C} we visualize the eigenvectors $v_j\in \C$, $j\neq 0$, in the complex plane for some $N\in \N$.
		\begin{figure}[!htb]
			\centering
			\subfloat[][\centering $N = 6$]{
				\label{fig:EV_C-a}
				\begin{minipage}[t]{.3\linewidth}
					\centering
					\includegraphics[width=0.9\textwidth]{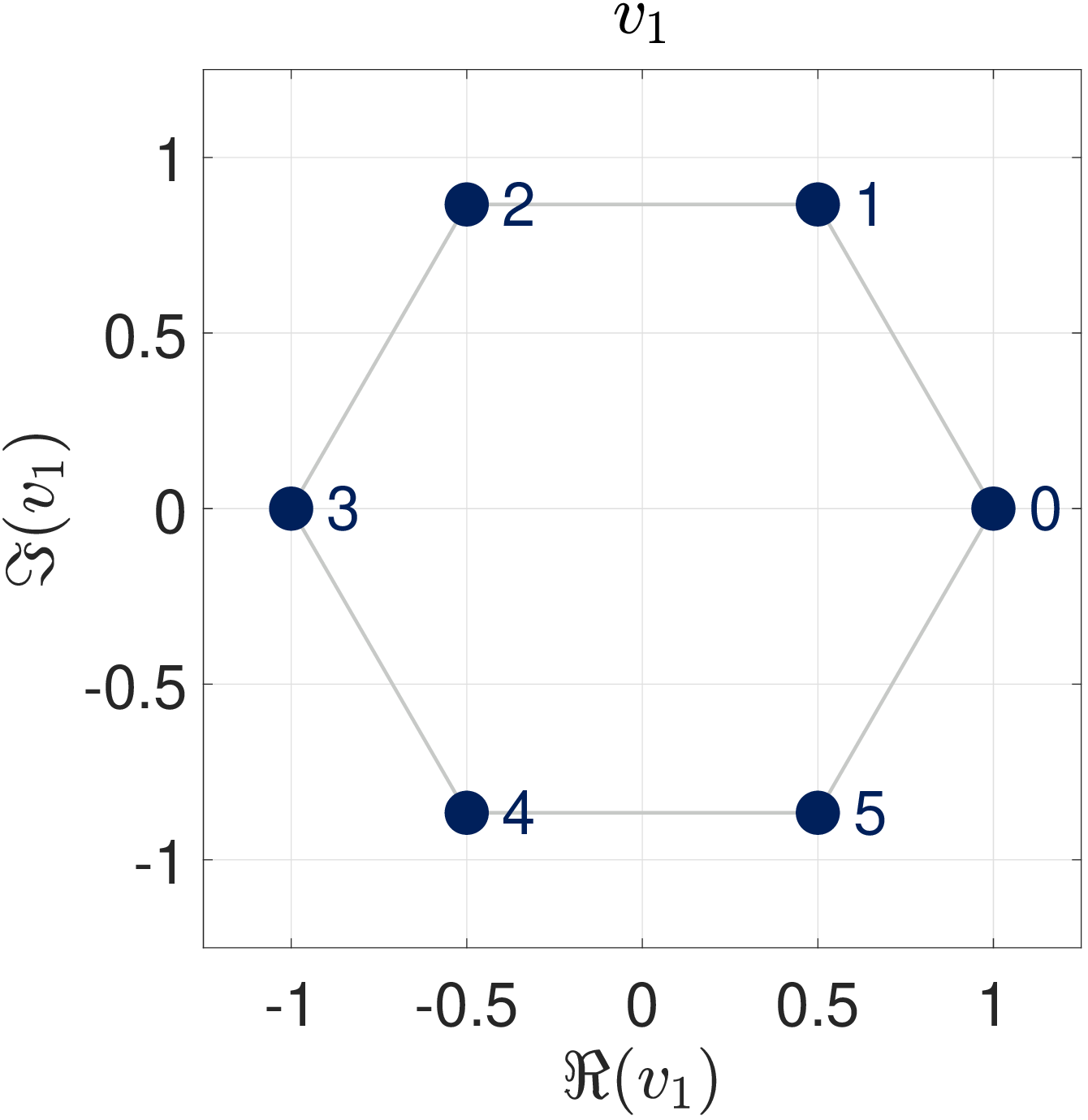}
				\end{minipage}
				\hfil
				\begin{minipage}[t]{.3\linewidth}
					\centering
					\includegraphics[width=0.9\textwidth]{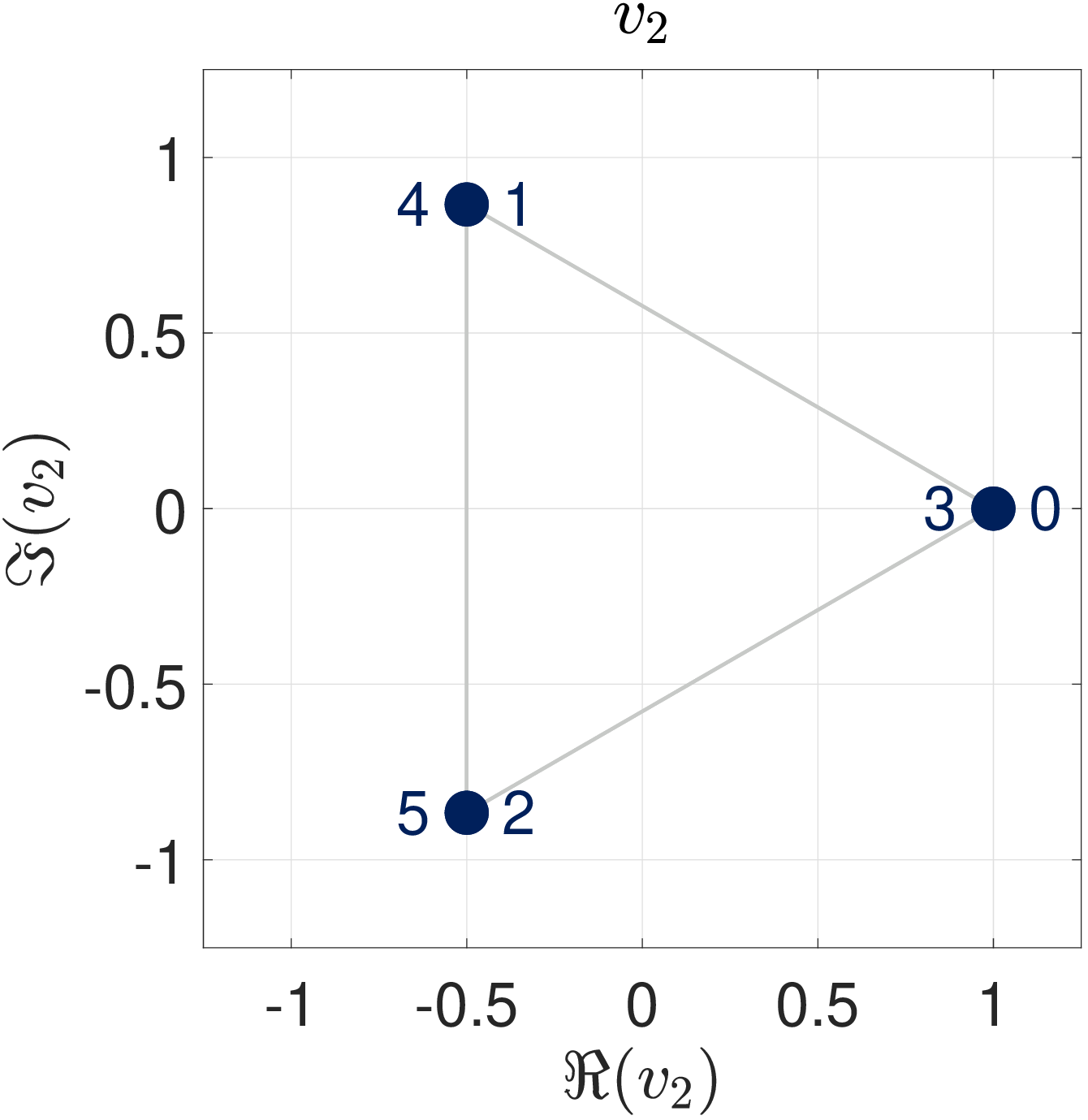}
				\end{minipage}	
				\hfil
				\begin{minipage}[t]{.3\linewidth}
					\centering
					\includegraphics[width=0.9\textwidth]{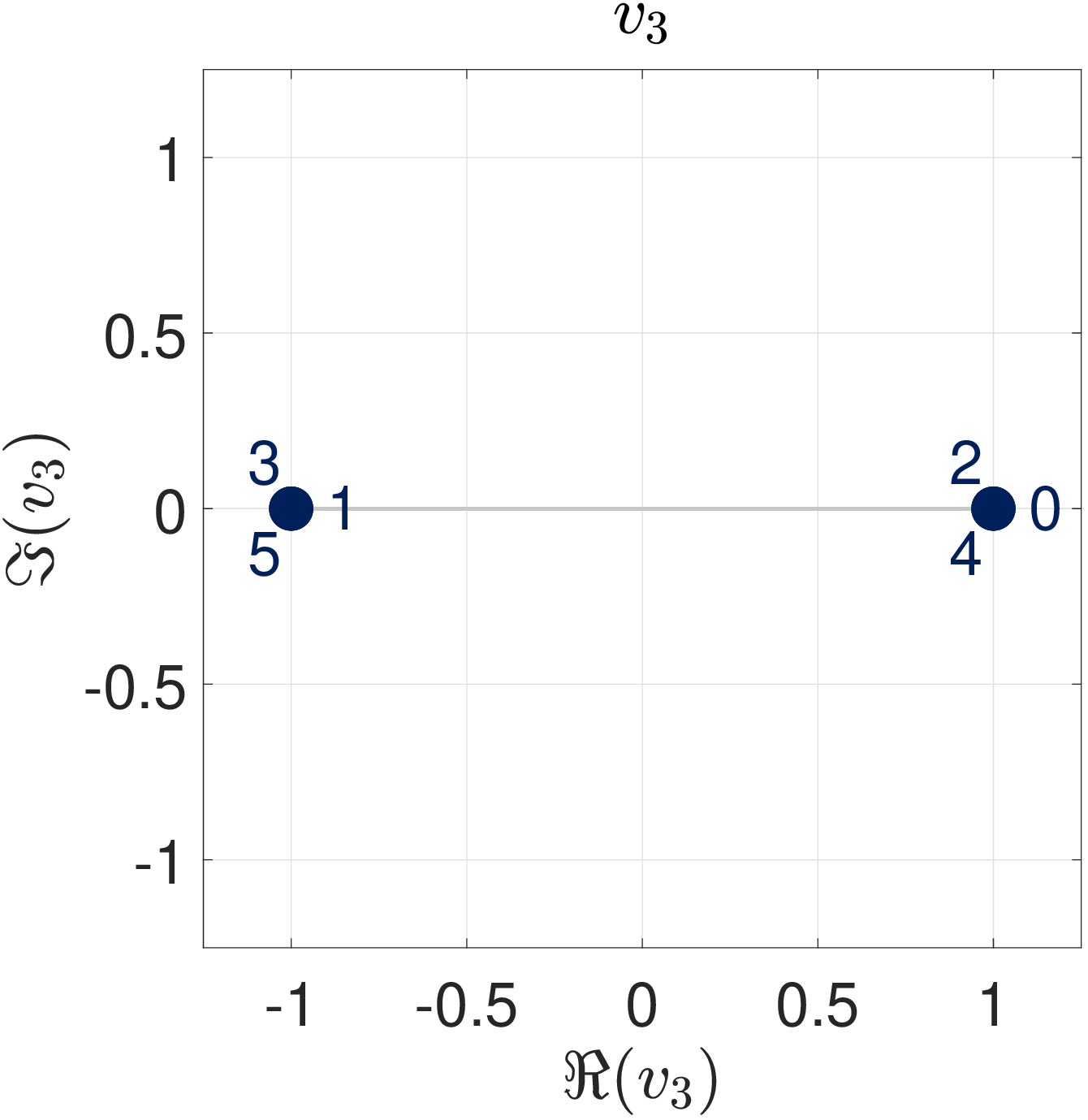}
			\end{minipage}			}\\
			\subfloat[][\centering $N = 7$]{
				\label{fig:EV_C-b}
				\begin{minipage}[t]{.3\linewidth}
					\centering
					\includegraphics[width=0.9\textwidth]{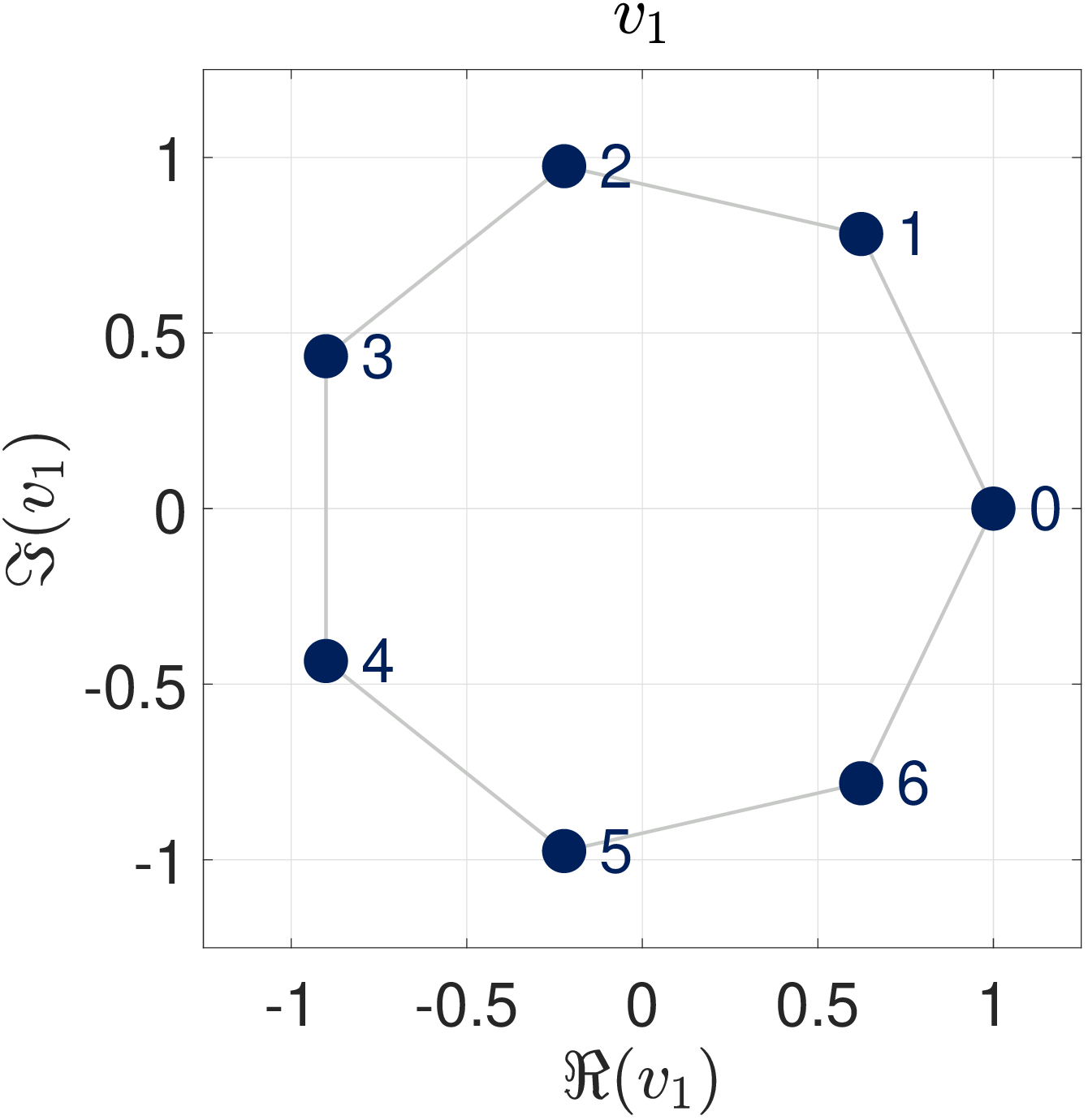}
				\end{minipage}
				\hfil
				\begin{minipage}[t]{.3\linewidth}
					\centering
					\includegraphics[width=0.9\textwidth]{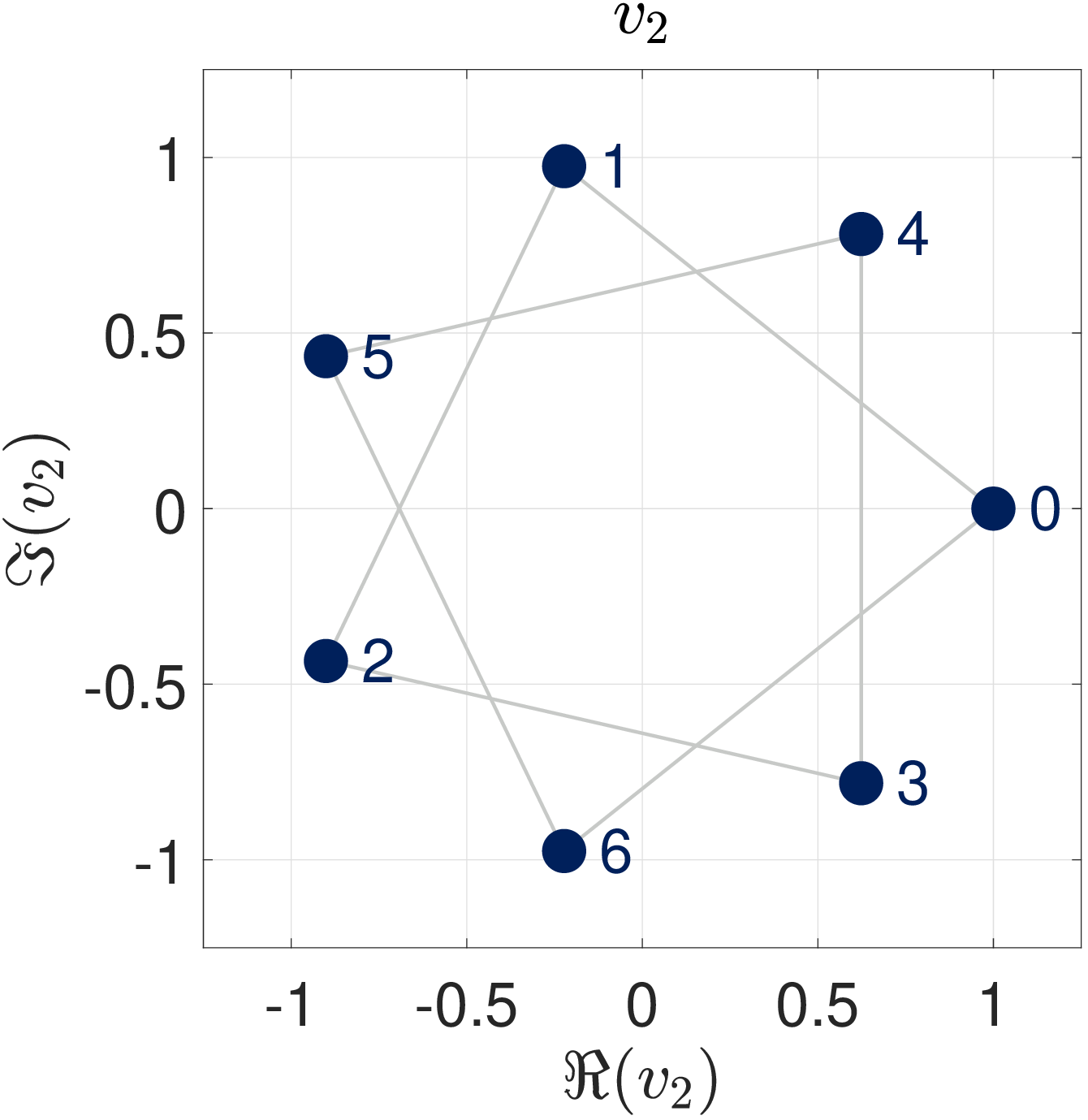}
				\end{minipage}	
				\hfil
				\begin{minipage}[t]{.3\linewidth}
					\centering
					\includegraphics[width=0.9\textwidth]{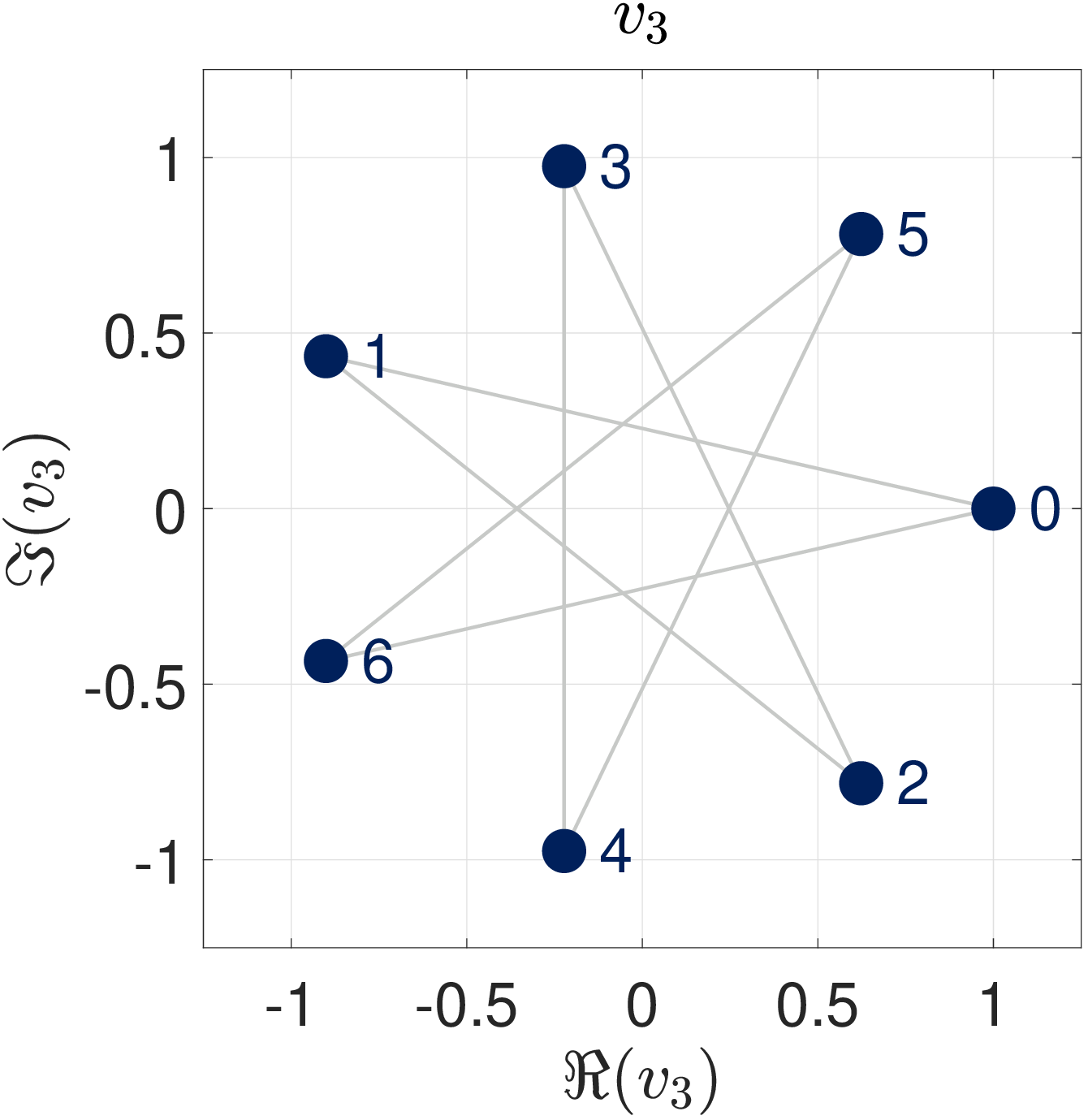}
				\end{minipage}			
			}    
			\caption{Visualization of the complex-valued eigenvectors $v_j\in \C^N$, $j\neq 0$, of a circulant matrix in the complex plane for (a) $N=6$ and (b) $N=7$. 
				Here, the number $i$ next to the blue dots correspond to $i$-th entry of $v_j\in \C$. The connections in gray between the points illustrate the (cyclically) ascending index $i = 0,1,\dots,N-1$. In particular, it does not indicate any dynamical influences according to an underlying interaction graph.
				Note that, by definition all points are on the unit circle and thus form a regular $N$-gon (see~\eqref{eq:EV}). For increasing $j$, if $N=2k$ is even, the $N$-gon 'degenerates' since some entries are the same, while for $N=2k+1$ only the ordering changes.
				Later on, we assign to every entry of $v_j\in \C$ its corresponding robot such that the dots also represent the position $x_i = \Re((v_j)_i)$ and $y_i = \Im((v_j)_i)$ of the $i$-th robot in the Euclidean plane (by changing the axis labels accordingly). Hence, the \emph{generating} configuration of the stable invariant subspace $V_j\subseteq \R^{2N}$ is also visualized
				(cf. \eqref{eq:Vj}).
			}
			\label{fig:EV_C}
		\end{figure}
	\end{enumerate}
\end{remark}

For the upcoming analysis, we derive solutions of $\eqref{eq:Z_lin}$ in closed form, which allows us to hierarchically rank them according to their speed of convergence.
By \Cref{prop:circ_spectral} a circulant matrix $W$ is diagonalizable, which yields the \emph{fundamental solutions} of \eqref{eq:Z_lin} to have the simple form
\begin{align*}
	{Z}_{j}^x(t) = e^{(-1+\lambda_j)t}\begin{pmatrix}v_{j}\\0\end{pmatrix} \text{ and }
	{Z}_{j}^y(t) = e^{(-1+\lambda_j)t}\begin{pmatrix}0\\v_{j}\end{pmatrix}, \quad j=0,\dotsc,N-1.
\end{align*}
However, for $j\neq 0$, the eigenvector $v_j$ and its corresponding eigenvalue $\lambda$ are, in general, complex-valued, which gives little possibility for interpretation in terms of configurations of $N\in \N$ robots
in the two-dimensional Euclidean plane.

We propose a real-valued basis as follows.
Let $N=2k$ be even or $N=2k+1$ be odd. For $j=0,\dotsc,k$ we define 
\begin{align}\label{eq:Vj}
	V_j = \Span\left(
	\begin{pmatrix}\Re(v_j)\\ -\Im(v_j)\end{pmatrix}, \begin{pmatrix}\Im(v_j)\\\Re(v_j)\end{pmatrix}, \begin{pmatrix}-\Im(v_j)\\ \Re(v_j)\end{pmatrix},\begin{pmatrix}\Re(v_j)\\ \Im(v_j)\end{pmatrix}
	\right),
\end{align}
as the subspaces
spanned by linear combinations of real- and imaginary parts of the eigenvector $v_j\in \C$. The collection of all $V_j\subseteq \R^{2N}$ decompose the state space as $\R^{2N} = \bigoplus_{j=0}^{k} V_j$ since the basis vectors are linearly independent.

The first three basis vectors in \eqref{eq:Vj} (considered as points/robots in the Euclidean plane) are reflections of the last one. In this sense, we say $v_j\in \C^N$, resp. the configuration $\left(\Re(v_j),\Im(v_j)\right)^T\in \R^{2N}$, \emph{generates} the subspace $V_j\subseteq \R^{2N}$. In particular, the $x$-coordinate of the $i$-th robot is given by the real part of the $i$-the entry of $v_j\in \C^N$, while its $y$-coordinate is the corresponding imaginary part.
Hence, the generating configuration $\left(\Re(v_j),\Im(v_j)\right)^T\in \R^{2N}$ of $V_j\subseteq \R^{2N}$ can be visualized as in \Cref{fig:EV_C} by changing the axis labels to $x$ and $y$.

On each $V_j\subseteq \R^{2N}$ the matrix $\mathbf{{W}}\in \R^{2N\times 2N}$ in \eqref{eq:Z_lin} acts as a block-diagonal matrix of two identical blocks $\Lambda_j\in \R^{2\times 2}$, i.e,
\begin{align}\label{eq:block_V}
	\mathbf{{W}}V_j = V_j\begin{pmatrix}\Lambda_j & \\ & \Lambda_j\end{pmatrix} \text{ with } \Lambda_j = \begin{pmatrix}\Re(\lambda_j)& \Im(\lambda_j)\\ -\Im(\lambda_j) &\Re(\lambda_j)\end{pmatrix}.
\end{align}
This is due to the fact, that the eigenvalue pair $\lambda_j \in \C$ and $\overline{\lambda_j} = \lambda_{N-j}$ act as a matrix of the form $\Lambda_j \in \R^{2\times 2}$. Here, we abuse notation and also consider $V_j$ to be the matrix containing the basis vectors of $V_j$ \eqref{eq:Vj}.
Since $\mathbf{W}\in \R^{2N\times 2N}$ maps each $V_j\subseteq \R^{2N}$ into itself, any initial configuration starting in one $V_j$ will remain therein for all times $t>0$. Moreover, the corresponding eigenvalues have real part less than $1$, which implies that solutions will furthermore converge to $0$ with the exponential decay rate $-1+\Re(\lambda_j)$ in $t$ (cf. \Cref{app:decomp}).
The \emph{smaller} $\Re(\lambda_j) < 1$ is, the \emph{faster} the solution converges. In particular, if $\Re(\lambda_s) < \Re(\lambda_j)$ for all $j\neq s$, the convergence rate of any configuration in $V_s\subseteq \R^{2N}$ is faster than any configuration in any other $V_j \subseteq \R^{2N}$ for $0\neq j\neq s$.

\begin{remark}
	If $v_j\in \C^N$ has only real-valued entries, a suitable basis of $V_j\subseteq \R^{2N}$ is already given by the first (or last) two spanning vectors and $\eqref{eq:block_V}$ reduces to $\mathbf{W}V_j = V_j\Lambda_j$. This is the case when (1) $j=0$ or (2) $N=2k$ and $j=k$. As it complicates the notation, we continue with the general case below.
	However, similar arguments can be made in this situation and we refer to 
	\Cref{app:decomp}
	for details.
\end{remark}

In order to dynamically rank the subspaces $V_j\subseteq \R^{2N}$ and obtain a dynamical hierarchy of gathering speeds, we define the following.
\begin{definition}\label{def:stable}
	Consider a gathering \cld/ protocol with weight matrix $W\in \R^{N\times N}$, eigenvalues $\lambda_j\in \C$ and invariant subspaces $V_j$ as in \eqref{eq:Vj}.
	\begin{enumerate}[label = (\alph*)]
		\item $\Syn \subseteq \R^{2N}$ denotes the $2$-dimensional subspace of all gathering points (cf. \Cref{def:gathering} and \Cref{rem:circ_EW} (a)).
		\item For $j\neq 0$, $V_j\subseteq \R^{2N}$ is called \emph{stable invariant subspace} with \emph{convergence rate} $\Re(\lambda_j) < 1$ (since the \cld/ protocol is gathering). 
		\item If $\Re(\lambda_s) < \Re(\lambda_j)$ for all $j\neq s$,  the subspace $V_s\subseteq \R^{2N}$ is called \emph{strong stable invariant subspace}.
	\end{enumerate}
\end{definition}

We summarize the results found above in the following theorem and state that the $2N$-dimensional state space of a gathering \cld/ protocol can be hierarchically decomposed. For more details of the proof, we refer to 
\Cref{app:decomp}.

\begin{theorem}\label{thm:decomp}
	For a gathering \cld/ protocol of $N\in \N$ robots with weight matrix $W\in \R^{N\times N}$ there exist a family of stable invariant subspaces $\left(V_j\right)_{j=1}^k\subseteq \R^{2N}$ with convergence rates $\Re(\lambda_j) < 1$ 
	and a $2$-dimensional subspace $\Syn\subseteq \R^{2N}$ of gathering points such that
	\begin{align}\label{eq:decomp}
		\R^{2N} = \Syn \oplus \left(\bigoplus_{j=1}^k V_j\right), \text{ where either}
	\end{align}
	\vspace{-1em}
	\begin{enumerate}[label = (\roman*)]
		\item  every subspace $V_j \subseteq \R^{2N}$ is $4$-dimensional, if $N=2k+1$ is odd,
		\item or,  if $N=2k$ is even, the subspace $V_k \subseteq \R^{2N}$ is $2$-dimensional, whereas the remaining $V_j\subseteq \R^{2N}$, $j\neq k$, are $4$-dimensional.
	\end{enumerate}
	Each subspace $V_j\subseteq \R^{2N}$ is spanned by a generating eigenvector $v_j \in \C^N$ of the weight matrix $W\in \R^{N\times N}$ in the sense of \eqref{eq:Vj}.
	By \Cref{rem:circ_EW} (c)
	this decomposition in \eqref{eq:decomp} is identical for every gathering \cld/ protocol. Only the explicit values of the convergence rates $\Re(\lambda_j)$ depend on the protocol itself.
\end{theorem}

\begin{example}\label{ex:EW}
	We apply \Cref{thm:decomp} to the running \Cref{ex:running_ex}.
	To this end, we apply the eigenvalue formula \eqref{eq:EW} to ${w=(0,1,0,\dotsc,0)^T \in \R^N}$ (cf. \Cref{ex:circ_gathering}):
	\begin{align}\label{eq:lambdaj}
		\lambda_j = \omega^j = \cos\left(\frac{2\pi j}{N}\right) + \mathbf{i}\sin\left(\frac{2\pi j}{N}\right)
	\end{align}
	and thus $\Re(\lambda_j) = \cos\left(\frac{2\pi j}{N}\right)$ for $j=0,1,\dotsc,N-1$.     These convergence rates are plotted in \Cref{fig:conv_rates} for some choices of $N\in \N$.
	\begin{figure}[!htb]
		\centering
		\subfloat[][\centering $N=6$]{
			\label{fig:conv_rates_d}
			\begin{minipage}[t]{.3\linewidth}
				\centering
				\includegraphics[width=0.9\textwidth]{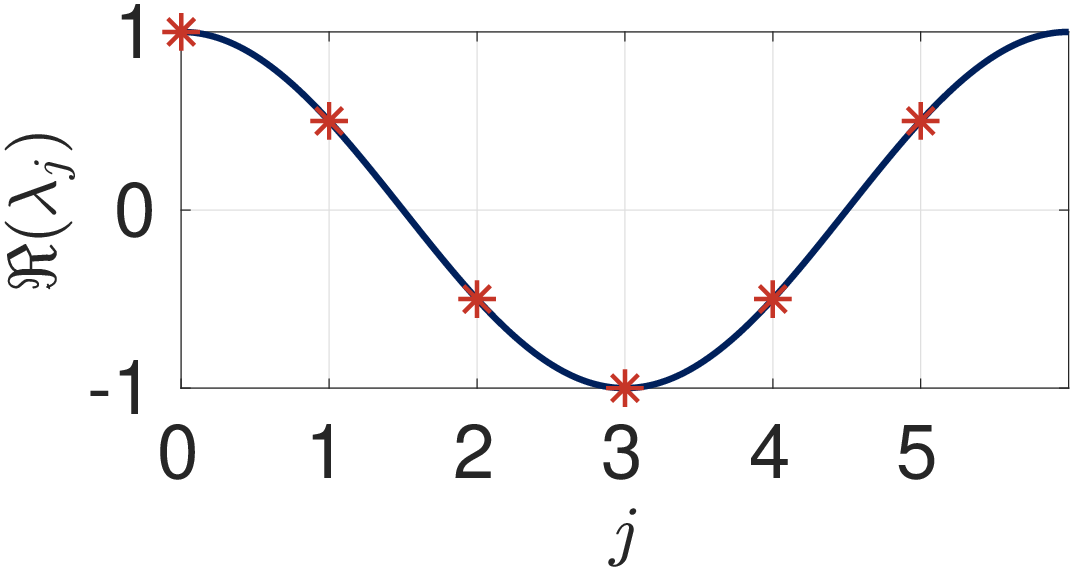}
			\end{minipage}
		}
		\hfil
		\subfloat[][\centering $N=7$]{
			\label{fig:conv_rates_e}
			\begin{minipage}[t]{.3\linewidth}
				\centering
				\includegraphics[width=0.9\textwidth]{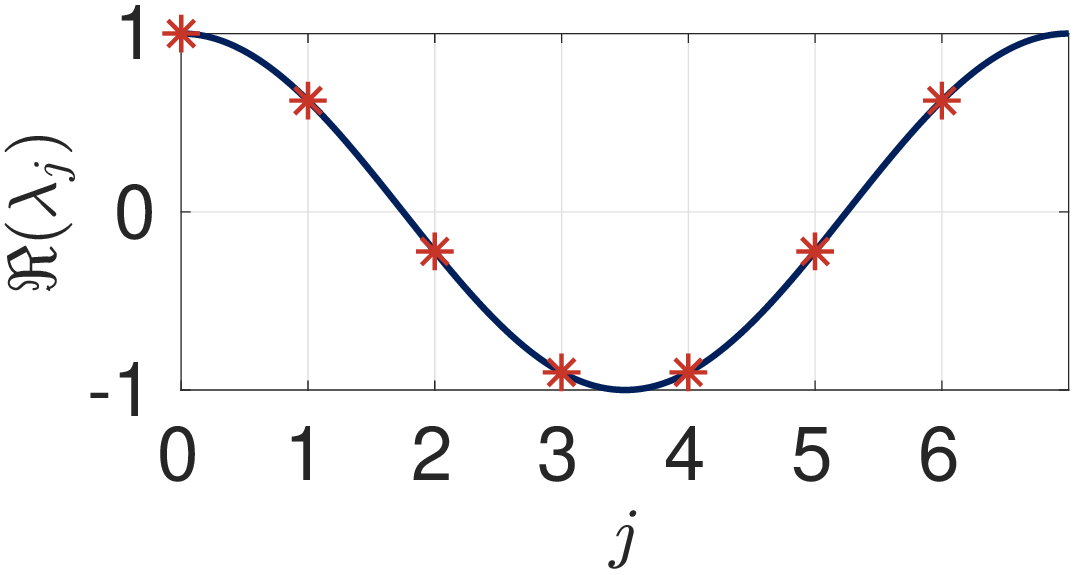}
			\end{minipage}
		}
		\hfil
		\subfloat[][\centering $N=8$]{
			\label{fig:conv_rates_f}
			\begin{minipage}[t]{.3\linewidth}
				\centering
				\includegraphics[width=0.9\textwidth]{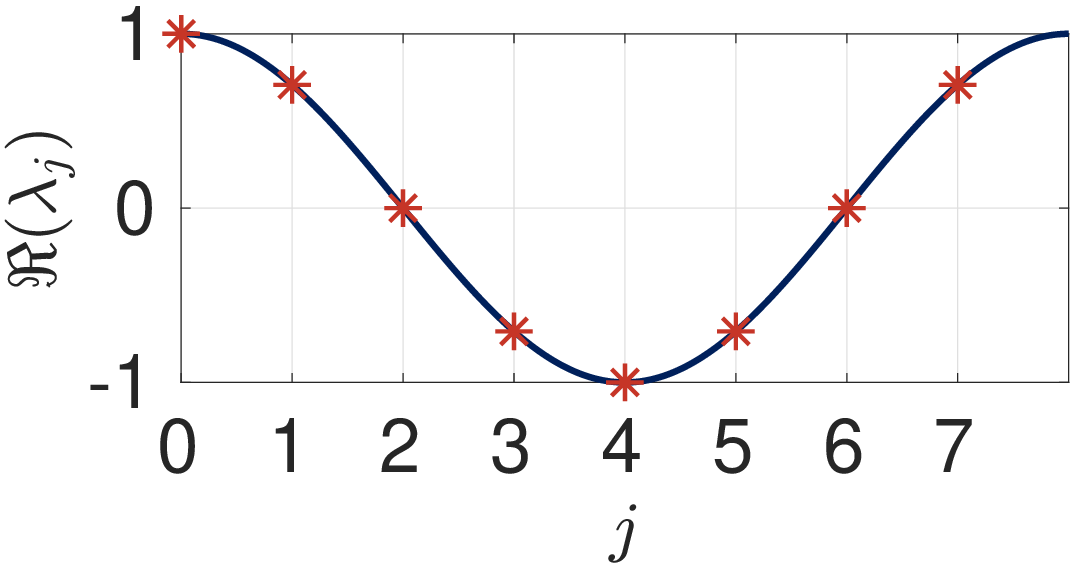}
			\end{minipage}
		}
		\caption{Illustration of the convergence rates $\Re(\lambda_j)$ of the \nbug{} problem for $N\in \set{6,7,8}$ by red stars. For reference, the cosine curve is plotted in blue. By \Cref{rem:circ_EW} (b) the convergence rates $\Re(\lambda_j)$ are symmetrically distributed.
		}
		\label{fig:conv_rates}
	\end{figure}
	
	Note that the convergence rates are strictly decreasing in $j$ (up to $k$), i.e., $\Re(\lambda_{j}) > \Re(\lambda_{j+1})$ for $j<k$. Hence, the stable subspaces $V_j\subseteq \R^{2N}$ ($j\neq 0$) are hierarchically ordered such that any configuration $Z \in V_j$ converges faster to a gathering point than any other configuration $\T{Z} \in V_i$ for $i<j$. 
	In particular, this model has a strong stable subspace $V_k \subseteq \R^{2N}$. 
	It is $4$-dimensional with convergence rate $\Re(\lambda_k) = \cos\left(\frac{2\pi k}{2k+1}\right)$ for $N=2k+1$ odd and $2$-dimensional with $\Re(\lambda_k) = -1$ for $N=2k$ even. The corresponding generating configuration is illustrated in the right panels of \Cref{fig:EV_C,fig:EV_C}.
\end{example}

As an arbitrary initial configuration ${Z}\in \R^{2N}$ can be written as a linear combination of the basis vectors of all $V_j\subseteq \R^{2N}$ and the \cld/ protocol is linear, \Cref{thm:decomp} implies that for $j\neq 0$ every part in $V_j\subseteq \R^{2N}$ vanishes as time proceeds and only the gathering point in $V_0 \subseteq \R^{2N}$ remains.

\begin{corollary}\label{cor:decomp}
	Let ${Z}(0) = (X(0),Y(0)) \in \R^{2N}$ be an initial configuration and the state space $\R^{2N}$ be decomposed as in \eqref{eq:decomp}. Then the solution ${Z}(t)\in \R^{2N}$ of the \cld/ protocol \eqref{eq:Z_lin} with initial condition ${Z}(0) \in \R^{2N}$ can be written as
	\begin{align}\label{eq:decomp_2}
		{Z}(t) = {Z}^* + \sum_{j=1}^k \alpha_j(t)~\Xi_j(t), \text{ where}
	\end{align}
	\begin{enumerate}[label = (\roman*)]
		\item ${Z}^* = (X^*,Y^*) = (x^*,\dotsc,x^*,y^*,\dotsc,y^*)$ is the final gathering point of ${Z} \in \R^{2N}$. In particular, ${x^* = \frac{1}{N} \sum_{i=0}^{N-1} X_i(0)}$ and  $y^* = \frac{1}{N} \sum_{i=0}^{N-1} Y_i(0)$.
		\item $\alpha_j(t) = \exp((-1+\Re(\lambda_j))t) \in \R$ is the exponentially decaying coefficient corresponding to $\Xi_j(t) \in V_j$.
		\item by abusing notation, $\Xi_j(t) = V_j \beta_j(t) \in V_j$ for a time-dependent 
		vector ${\beta_j(t) \in \R^{\dim V_j}}$ with constant norm, i.e., $\norm{\beta_j(t)}_2 = \norm{\beta_j(0)}_2$ for all $t \geq 0$.
	\end{enumerate}
\end{corollary}

\begin{example}
	We apply \Cref{cor:decomp} to the \nbug{} problem.
	For $N=7$, we show in \Cref{fig:init} a random initial configuration ${Z}(0) \in \R^{2N}$ as well as its final gathering point ${Z}^*\in \R^{2N}$. The corresponding exponentially decaying coefficients $\alpha_j(t) \in \R$ are plotted in \Cref{fig:alpha}.
	The initial decomposition into $\Xi_j(0)\in V_j$ is visualized in \Cref{fig:decomp-xi}. The coefficients $\beta_j(t)\in \R^{\dim V_j}$ are shown in \Cref{fig:decomp-beta}.  
	\begin{figure}[!htb]
		\centering
		\subfloat[][\centering ${Z}(0)\in \R^{2N}$ and ${Z}^*\in \R^{2N}$]{
			\label{fig:init}
			\begin{minipage}[t]{.45\linewidth}
				\centering
				\includegraphics[width=0.6\textwidth]{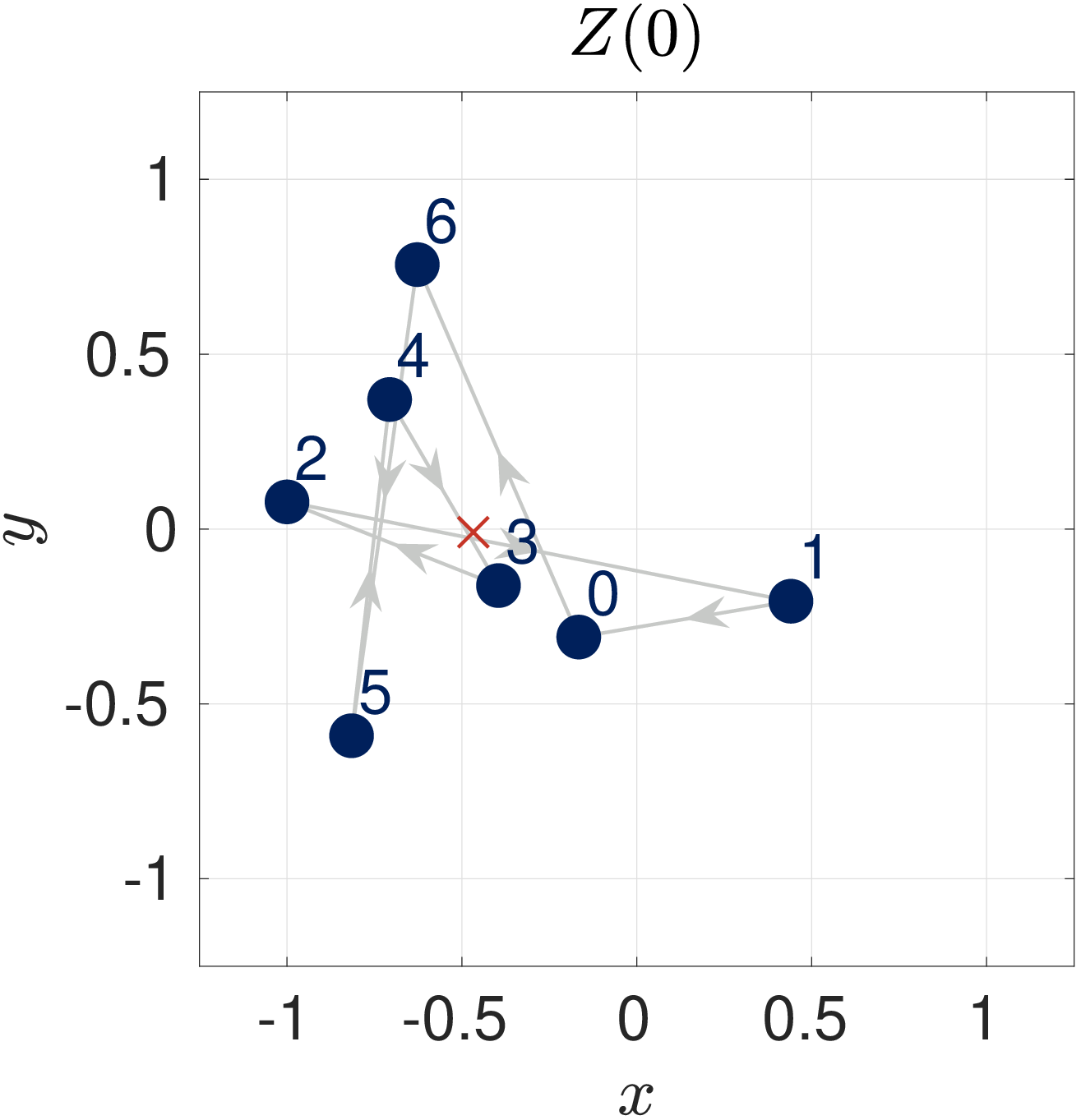}
			\end{minipage}
		}
		\hfil
		\subfloat[][\centering  Trajectory of $\alpha_j(t) \in \R$]{
			\label{fig:alpha}
			\begin{minipage}[t]{.45\linewidth}
				\centering
				\includegraphics[width=0.6\textwidth]{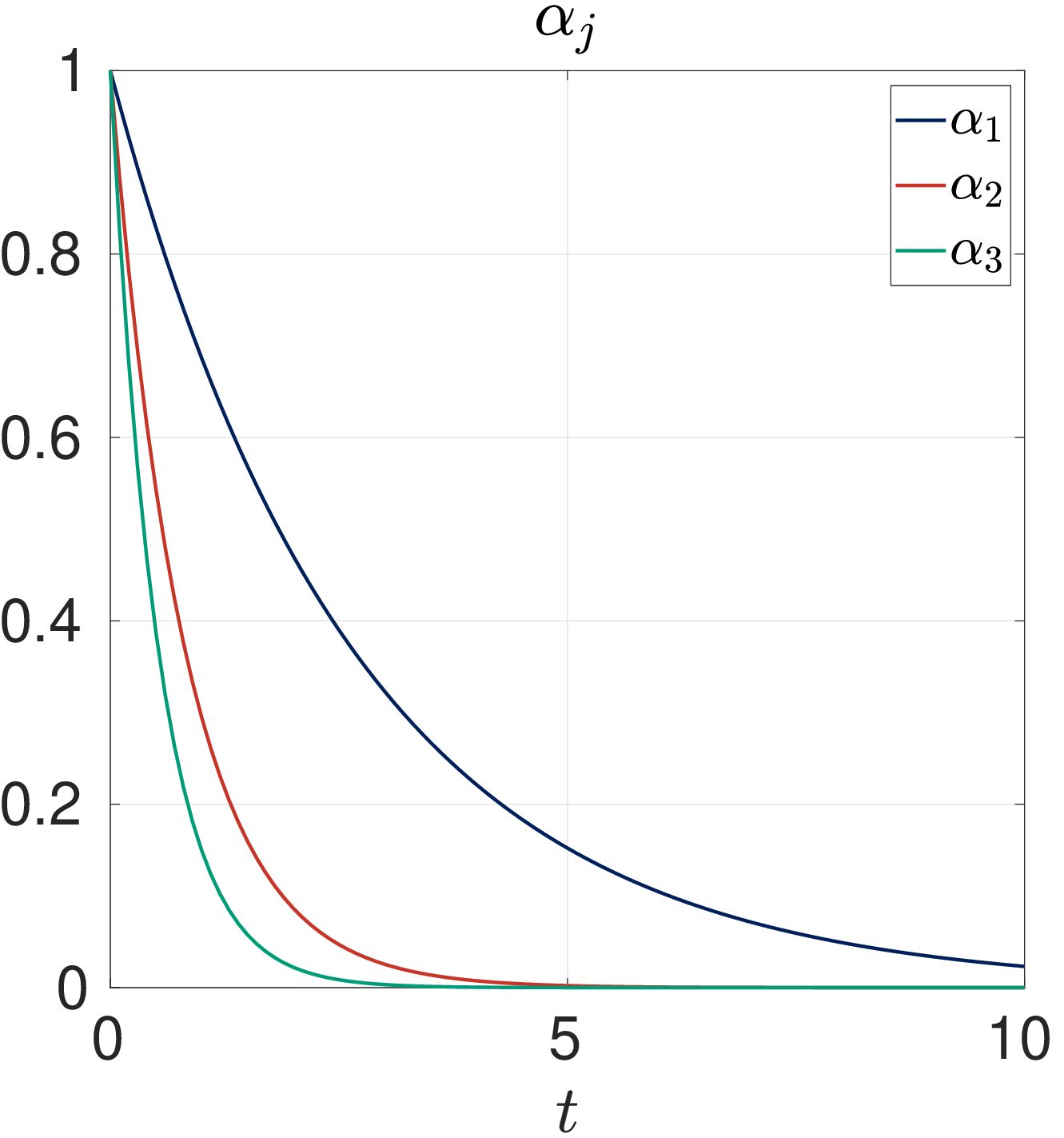}
			\end{minipage}
		}  
		\caption{Visualization of a random initial configuration ${Z}(0)\in \R^{2N}$ (blue dots) and its gathering point $Z^*\in \R^{2N}$ (red cross) in the Euclidean plane for $N=7$ in (a). Its corresponding decaying coefficients $\alpha_j(t) \in \R$ are shown in (b). Note that $\alpha_j(t) \in \R$ decreases faster for increasing $j$, which illustrates the dynamical hierarchy discussed in \Cref{ex:EW}.}
	\end{figure}
	\begin{figure}[!htb]
		\centering
		\subfloat[][\centering Initial decomposition of ${Z}(0)$ into $\Xi_j(0) \in V_j$]{
			\label{fig:decomp-xi}
			\begin{minipage}[t]{.3\linewidth}
				\centering
				\includegraphics[width=0.9\textwidth]{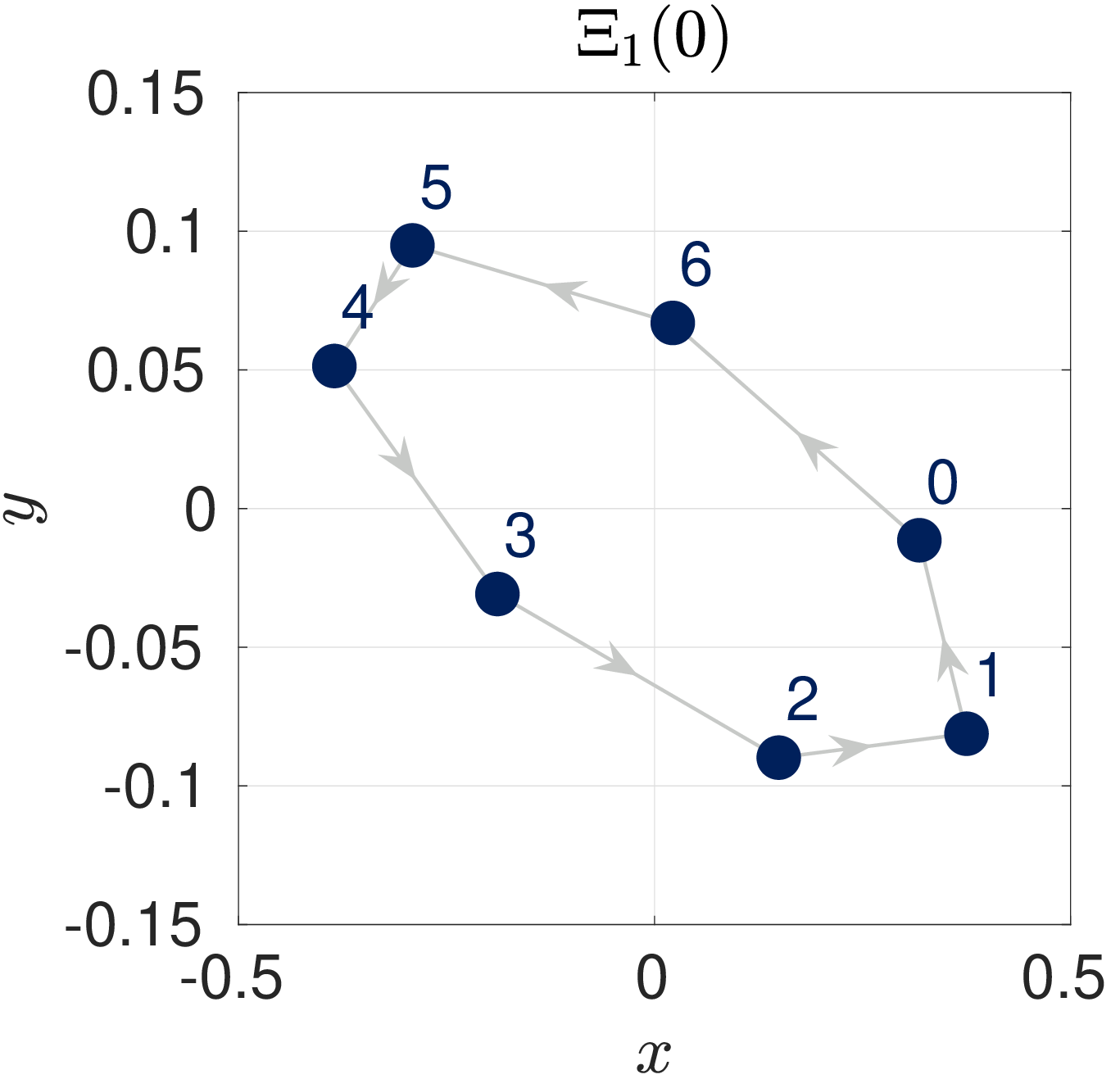}
			\end{minipage}
			\hfil
			\begin{minipage}[t]{.3\linewidth}
				\centering
				\includegraphics[width=0.9\textwidth]{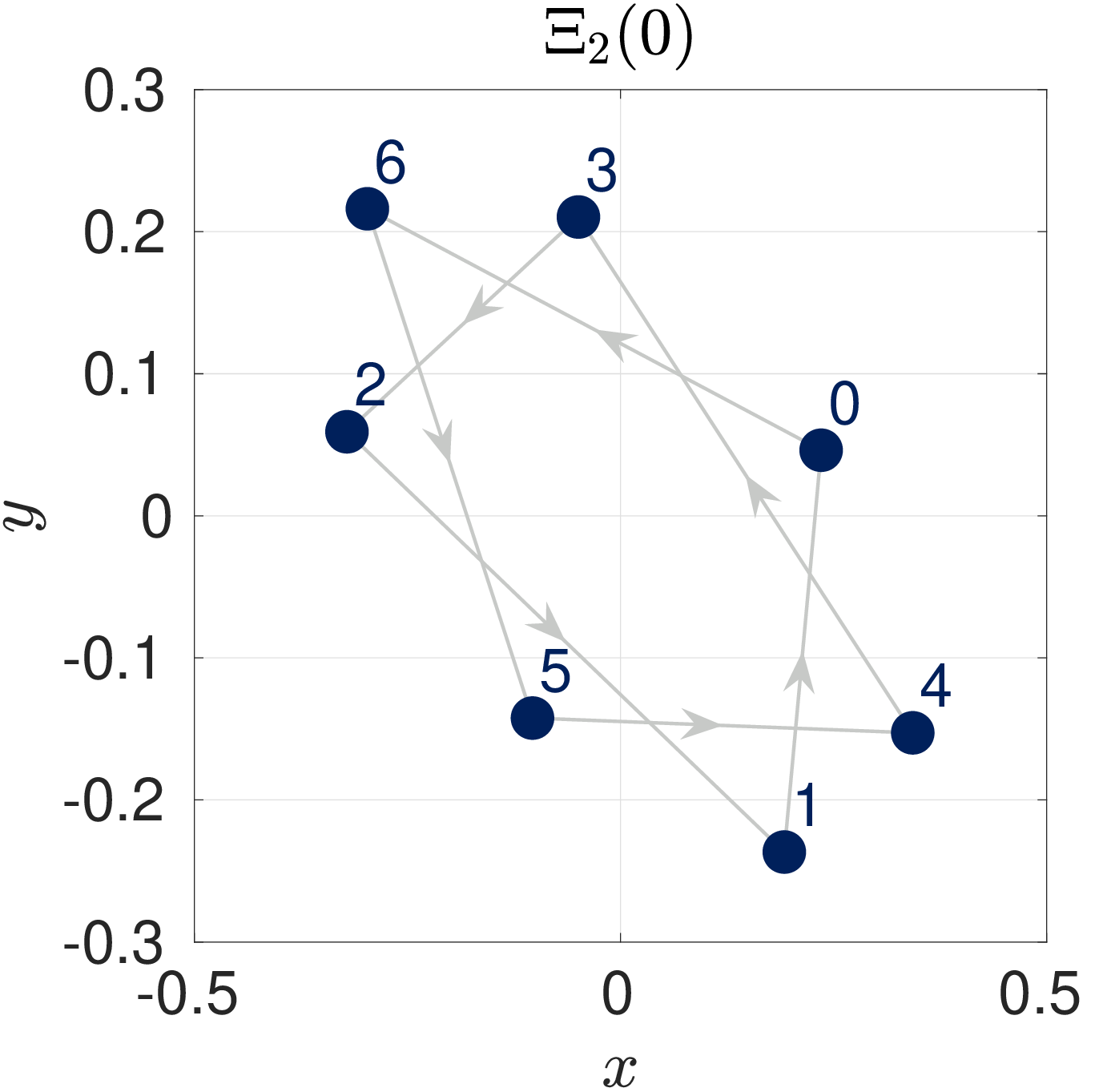}
			\end{minipage}	
			\hfil
			\begin{minipage}[t]{.3\linewidth}
				\centering
				\includegraphics[width=0.9\textwidth]{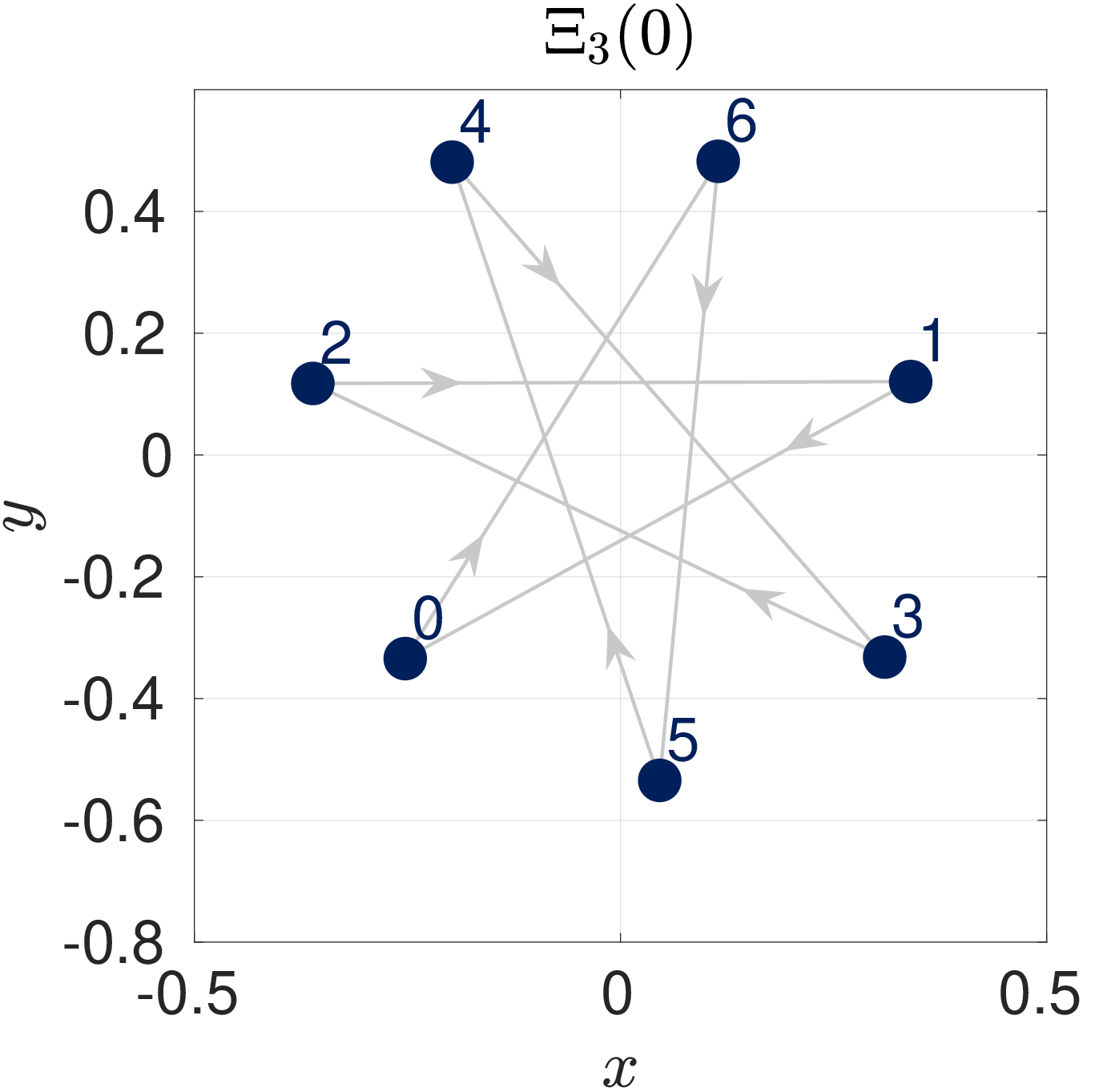}
			\end{minipage}
		}\\
		\subfloat[][\centering
		Trajectory of 
		$\beta_j(t) \in \R^{\dim V_j}$ such that $\Xi_j(t) = V_j \beta_j(t) \in V_j$
		]{
			\label{fig:decomp-beta}
			\begin{minipage}[t]{.3\linewidth}
				\centering
				\includegraphics[width=0.9\textwidth]{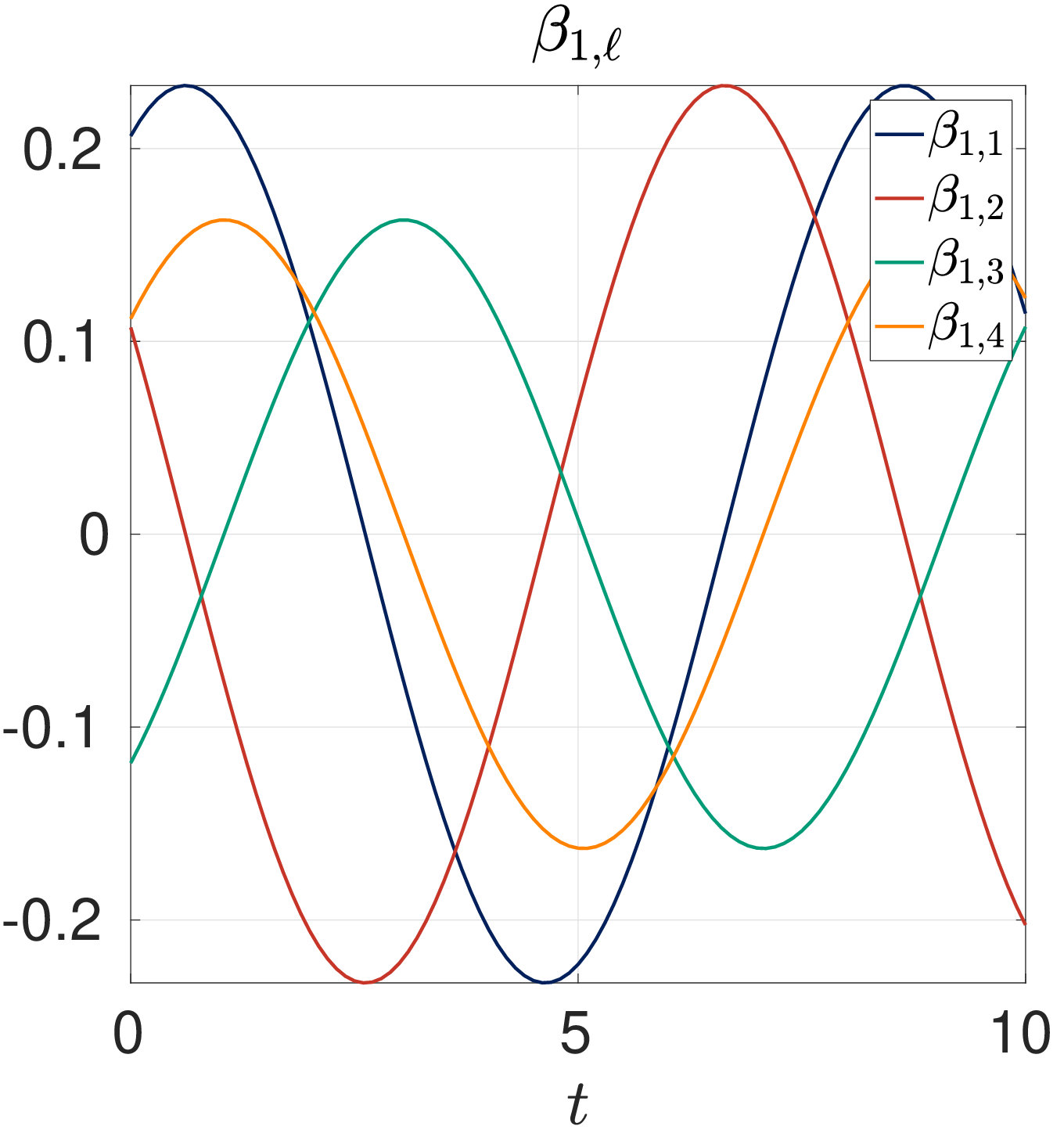}
			\end{minipage}
			\hfil
			\begin{minipage}[t]{.3\linewidth}
				\centering
				\includegraphics[width=0.9\textwidth]{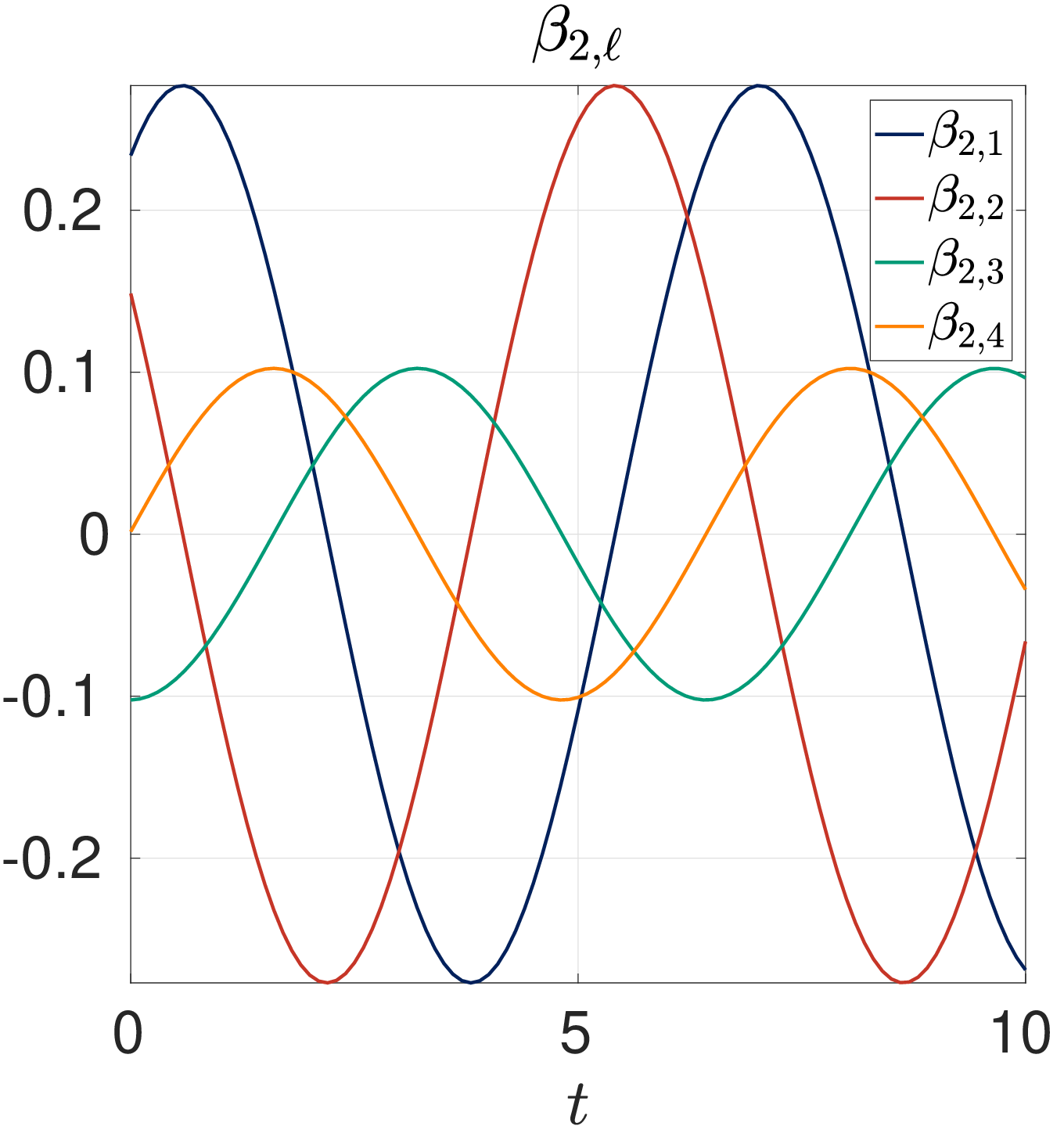}
			\end{minipage}
			\hfil
			\begin{minipage}[t]{.3\linewidth}
				\centering
				\includegraphics[width=0.9\textwidth]{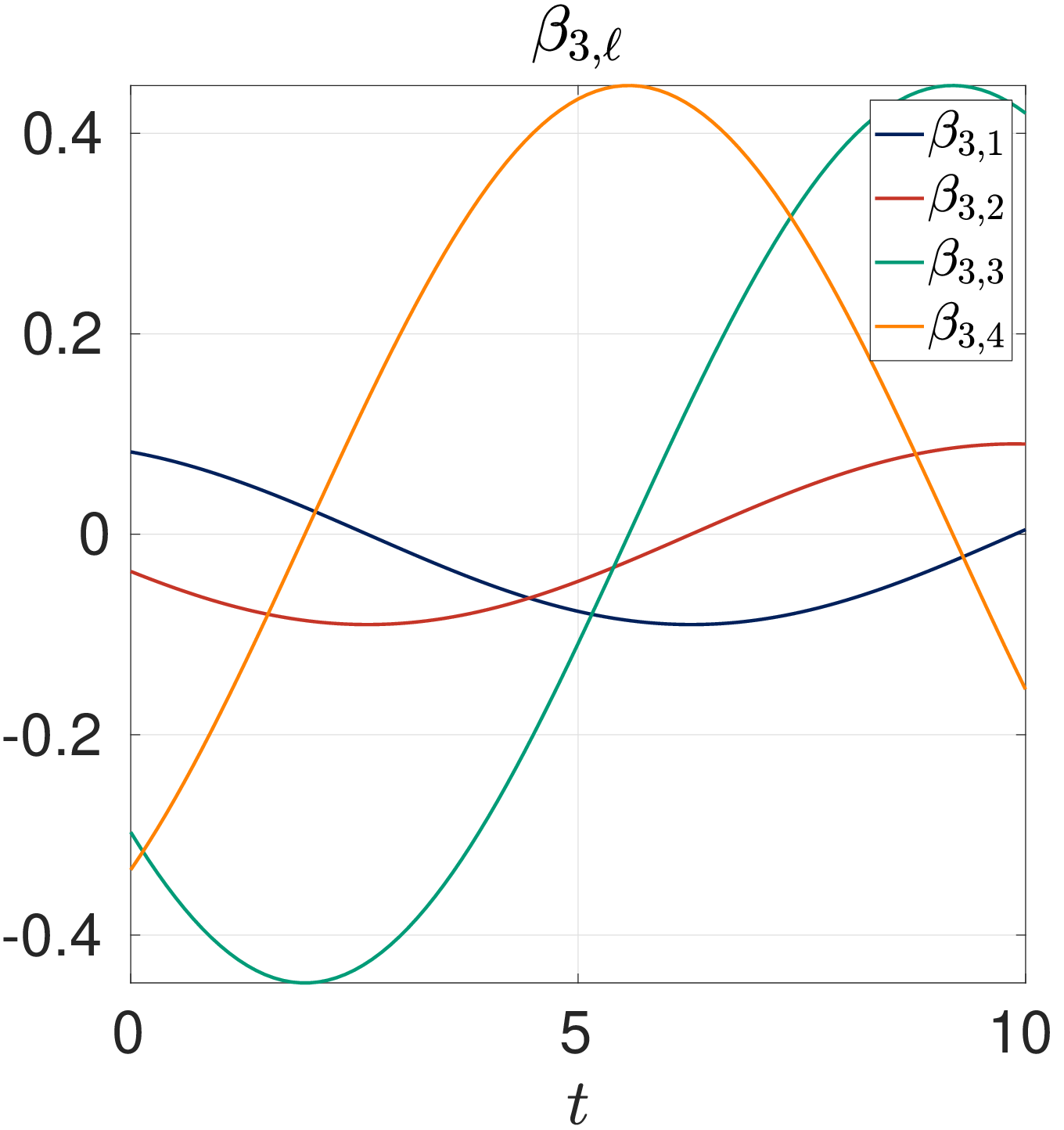}
			\end{minipage}
		} 
		\caption{Illustration of the decomposition \eqref{eq:decomp_2} and its individual dynamics. Note that $\beta_j(t)\in \R^{\dim V_j}$ only consists of sinusoidal curves since the decaying part is contained in $\alpha_j(t)\in \R$. In particular, we have $\norm{\beta_j(t)}_2 =\norm{\beta_j(0)}_2$ for all $t\geq 0$.}
		\label{fig:decomp}
	\end{figure}
\end{example}

\section{Analysis of the Nonlinear Scaled Circulant Laplacian Dynamics}
\label{sec:analysis_nonlin}
In this section, we study the dynamical properties of \nscld/ protocols of the form \eqref{eq:f_lin_normalized} by extending \Cref{thm:decomp} to the nonlinear case. Reconsider the subspaces $V_j \subseteq \R^{2N}$ for $j\neq 0$ that are invariant for \cld/ protocols.
Due to their block structure in \eqref{eq:block_V}, these can be further decomposed into $V_j = V_j^+ \oplus V_j^-$ with
\begin{align*}
	V_j^+ = \Span\left(
	\begin{pmatrix}\Re(v_j)\\ -\Im(v_j)\end{pmatrix},
	\begin{pmatrix}\Im(v_j)\\\Re(v_j)\end{pmatrix}\right)
	\text{ and } 
	V_j^- = \Span\left(
	\begin{pmatrix}-\Im(v_j)\\ \Re(v_j)\end{pmatrix},
	\begin{pmatrix}\Re(v_j)\\ \Im(v_j)\end{pmatrix}
	\right)
\end{align*}
on which we have $WV_j^\pm = V_j^\pm \Lambda_j$, where $V_j^\pm$ denotes either
$V_j^+$ or $V_j^-$.

We evaluate the right-hand side of \eqref{eq:f_lin_normalized} for ${Z} \in V_j^\pm$.
As computing $\NN$, resp. $\n$, in \eqref{eq:f_lin_normalized} involves evaluating the 'individual' length of each direction $(x_i,y_i)\in \R^2$ we define for ${Z} = (X,Y) \in \R^{2N}$ the length of the $i$-th direction by
\begin{align}\label{eq:normi}
	\normi({Z}) = \norm{(x_i,y_i)}_2.
\end{align}	
A short computation shows that $\normi$ has a useful representation when writing ${Z}\in V_j^\pm$ as a linear combination of the basis vectors
(see \Cref{app:decomp_nonlinear}).
Since the entries of $v_j$ are roots of unity we have
\begin{align*}
	{Z} = V_j^\pm \begin{pmatrix}a\\b\end{pmatrix} = a\begin{pmatrix}\Re(v_j)\\ \mp\Im(v_j)\end{pmatrix} + b \begin{pmatrix}\pm\Im(v_j)\\\Re(v_j)\end{pmatrix} & \implies \normi({Z}) = \norm{(a,b)}_2,
\end{align*}
which only depends on the coefficients $(a,b)\in \R^2$ and, in particular, \emph{not} on the index $i$. What remains to consider is how multiplication of $Z\in V_j^\pm$ by $-\mathbf{I}_{2N} + \mathbf{W}$ affects $\normi$. To this end, we compute
\begin{align*}
	\normi\left(\left(-\mathbf{I}_{2N}+\mathbf{\mathbf{W}}\right){Z}\right) = \norm{\left(-\mathbf{I}_2 + \Lambda_j\right) \begin{pmatrix}a\\b\end{pmatrix}}_2 = \abs{-1 + \lambda_j}\norm{(a,b)}_2,
\end{align*}
since $V_j^\pm \subseteq \R^{2N}$ is invariant under the linear dynamics and $\Lambda_j \in \R^{2\times 2}$ acts on $\R^2$ as the multiplication by $\lambda_j\in \C$. Again, this expression only depends on the eigenvalue $\lambda_j \in \C$ and the coefficients $(a,b)\in \R^2$.

As the length of each individual direction is the same for ${Z} \in V_j^\pm$, we can rewrite the scaling in \eqref{eq:f_lin_normalized} as the multiplication with a single non-negative scalar $\n(\abs{-1+\lambda_j}\norm{(a,b)}_2)\in \R$.
This yields that the dynamics of a \nscld/ protocol leaves the subspaces $V_j^\pm$ invariant and \eqref{eq:f_lin_normalized} can be written as
\begin{align}\label{eq:DS_TZ}
	\dot{{Z}} = \n(\abs{-1 + \lambda_j}\norm{(a,b)}_2) \left(-\mathbf{I}_{2N} + \mathbf{W}\right){Z}, \quad Z\in V_J^\pm.
\end{align}
Thus, the nonlinear dynamics restricted to $V_j^\pm \subseteq \R^{2N}$ reduces to a nonlinear uniform scaling of the linear dynamics.
In $(a,b)$-coordinates, this becomes
\begin{align}\label{eq:ab}
	\dot{\begin{pmatrix}a\\b\end{pmatrix}} = \n(\abs{-1 + \lambda_j}\norm{(a,b)}_2) \left(-\mathbf{I}_2 + \Lambda_j\right) \begin{pmatrix}a\\b\end{pmatrix}.
\end{align}
To derive an analogue of the convergence rate for \nscld/ protocols, we consider the dynamical behavior of $\norm{(a,b)}_2$ on $V_j^\pm$ for $j\neq 0$. Using the chain rule in combination with \eqref{eq:ab} its dynamics is given by
\begin{align}\label{eq:dotnorm}
	\dot{\norm{\begin{pmatrix}a\\b\end{pmatrix}}}_2 = (-1+\Re(\lambda_j))  \n(\abs{-1+ \lambda_j}\norm{(a,b)}) \norm{\begin{pmatrix}a\\b\end{pmatrix}}_2 < 0,
\end{align}
since $\Re(\lambda_j)< 1$ for $j\neq 0$ (cf. \Cref{thm:decomp}) and $\n$ is assumed to be non-negative
(cf. \Cref{app:decomp_nonlinear} for details).
We summarize our findings in the following theorem.
\begin{theorem}\label{thm:decomp_nonlinear}
	For an \nscld/ protocol with underlying gathering \cld/ protocol of $N\in \R$ robots with weight matrix $W\in \R^{N\times N}$ there exist two families of stable invariant subspaces $\left(V_j^+\right)_{j=1}^k\subseteq \R^{2N}$ and $\left(V_j^-\right)_{j=1}^k\subseteq \R^{2N}$ assigned with convergence rate $\Re(\lambda_j) < 1$, as well as is a $2$-dimensional subspace $\Syn\subseteq \R^{2N}$ of gathering points such that	
	\begin{align}\label{eq:decomp_nonlinar}
		\R^{2N} = \Syn \oplus \left(\bigoplus_{j=1}^k V_j^+ \oplus V_j^-\right), \text{ where either}
	\end{align}
	\begin{enumerate}[label = (\roman*)]
		\item every subspace $V_j^\pm \subseteq \R^{2N}$ is $2$-dimensional, if $N=2k+1$ is odd,
		\item or, if $N=2k$ is even, the subspace $V_k^\pm \subseteq \R^{2N}$ is $1$-dimensional, whereas the remaining $V_j^\pm\subseteq \R^{2N}$, $j\neq k$, are $2$-dimensional.
	\end{enumerate}
	In each of the spaces $V_j^\pm$, gathering properties are described by \eqref{eq:dotnorm}. Again, this decomposition is identical for every \nscld/ protocol.
\end{theorem}

In addition to that, it can be shown that visibility is preserved in $V_j^\pm$. More precisely, \Cref{prop:visibility-circulant-nonlinear} below states that interacting robots do not lose sight of each other, if the initial configuration was in any of the subspaces $V_j^\pm$ and all robots had all their interaction partners within viewing range. The statement is a consequence of the specific equation of motion on $V_j^\pm \subseteq \R^{2N}$ given by \eqref{eq:DS_TZ} and an analogous estimates as in the proof of \Cref{thr:visibility-circulant}
(details in \Cref{app:vis_nonlinear}).

\begin{proposition}
	\label{prop:visibility-circulant-nonlinear}
	Consider an \nscld/ protocol with underlying gathering \cld/ protocol of $N\in \R$ robots with weight matrix $W\in \R^{N\times N}$ and decomposition $\eqref{eq:decomp_nonlinar}$. 
	Let ${z_0(0),\dotsc,z_{N-1}(0)\in\R^2}$ be a valid initial configuration in one of the subspaces $V_j^\pm \subseteq \R^{2N}$, that is,
	$\|z_p(0)-z_q(0)\| \le \mathcal{C}$ for all $(q,p)\in E$.
	Then, visibility is preserved under the dynamics, i.e.,
	\[ \|z_p(t)-z_q(t)\| \le \mathcal{C} \quad\text{for all } t\geq 0 \text{ and } (q,p)\in E. \]
\end{proposition}

Finally, we discuss the implication of \Cref{thm:decomp_nonlinear} for the non-smooth scaling
$\n(\norm{(x,y)}_2 = \norm{(x,y)}_2^{-1}$
in comparison to no scaling, i.e, the linear case.
Non-smooth scaling causes $Z\in V_j^\pm \subseteq \R^{2N}$ to gather in \emph{finite time}
depending on the eigenvalue $\lambda_j \in \C$  and the 'size' of $Z$.
For the proof we refer to 
\Cref{app:finite_time}.

\begin{proposition}\label{prop:finite_time}
	An \nscld/ protocol with $\n(\norm{(x,y)}_2) = \norm{(x,y)}_2^{-1}$ achieves gathering of initial configurations $Z_0\in V_j^\pm \subseteq \R^{2N}$ ($j\neq 0$) in finite time that, for any $i$, is given by
	\begin{align*}
		t_{\mbox{gather},j} = -\frac{\abs{-1+\lambda_{j}}}{-1+\Re(\lambda_j)} \normi(Z_0) = -\frac{\abs{-1+\lambda_{j}}}{-1+\Re(\lambda_j)}\norm{\begin{pmatrix}a_0\\b_0\end{pmatrix}}_2 > 0.
	\end{align*}
\end{proposition}

In contrast, without scaling, i.e, purely linear dynamics, \eqref{eq:dotnorm} becomes
\begin{align*}
	\dot{\norm{\begin{pmatrix}a\\b\end{pmatrix}}}_2  = (-1+\Re(\lambda_j))\norm{\begin{pmatrix}a\\b\end{pmatrix}}_2.
\end{align*}
Convergence is exponential with decay rate $-1+\Re(\lambda_j)$. However, the origin is never actually reached, meaning that gathering is not possible in finite time.

\begin{example}
	According to \eqref{eq:lambdaj} and \Cref{prop:finite_time} the gathering times for the \nbug{} problem with scaling $\n(\norm{(x,y)}_2) = \norm{(x,y)}_2^{-1}$ are given by
	\begin{align*}
		t_{\mbox{gather},j} =  -\frac{\abs{-1+\lambda_{j}}}{-1+\Re(\lambda_j)} = -\frac{\abs{-1 +\cos\left(\frac{2\pi j}{N}\right) + \mathbf{i}\sin\left(\frac{2\pi j}{N}\right)}}{-1+\cos\left(\frac{2\pi j}{N}\right)},
	\end{align*}
	when a normalized initial condition $Z\in V_j^\pm$ is considered, i.e., $\norm{(a_0,b_0)}_2 = 1$. For $N=2k$ we have the real eigenvalue $\lambda_k = -1$ corresponding to the \emph{strong} stable subspace $V_k^\pm \subseteq \R^{2N}$ 
	and the \emph{fastest} gathering time is $t_{\mbox{gather},k} = 1$.
\end{example}

\section{Conclusion and Outlook}
\label{sec:conclusion}
In this article, we investigate gathering properties of \cld/ and \nscld/ protocols and derived a decomposition of the swarm state spaces into stable invariant subspaces which represent a dynamical hierarchy in terms convergence speed. The decomposition is independent of the explicit protocol, whereas only the gathering rates depend on the weights in the interaction graph.

We illustrate our findings on the \nbug{} problem which is in the class of our considered  protocols. Certainly, there are more protocols that can be described by \cld/ and/or \nscld/ models such as the \gtm{} and the \gta{}
protocol. The applications of our findings to these protocols can be found in 
\Cref{app:example}.

We consider fixed \emph{circulant} interaction structure. However, our results can be adapted to \emph{general} fixed interaction graphs. Here, the eigendecomposition is still crucial for the dynamical analysis and stable invariant subspaces given by eigenspaces with different gathering rates can be found. Although, the explicit structure of this decomposition is generally unclear as the spectral decomposition of an arbitrary weight matrix cannot be analytically derived.

In future research, the general case could be considered, where the interaction structure is not fixed, but rather depends on the \emph{visibility graph}, i.e, $e=(j,i)\in E$ if and only if robots $i$ and $j$ are sufficiently close. Here, a well-known model is the \oblot/ model~\cite{DBLP:series/lncs/FlocchiniPS19}, where robots are \emph{anonymous} in contrast to our \cmas/ model. In this situation the interaction graph is only changing finitely often at discrete times and we expect that the eigendecomposition (in between each topology change) remains a useful tool for the dynamical analysis.

Finally, in this work we consider a time-continuous robot model. However, a similar dynamical analysis is possible for a time-discrete model as well. In particular, the mathematical theory of linear discrete dynamical systems is just as well developed. In this context, it is an interesting open question if the results can even be extended to \emph{asynchronous} time models.

\section*{Acknowledgements}
The authors would like to express their special gratitude to
Friedhelm Meyer auf der Heide, Jannik Castenow, and Jonas
Harbig who provided invaluable inspiration and counsel for
the distributed computing side of the project and took part in
the development of some of the research ideas presented in this
manuscript in numerous spirited and productive discussions.

\bibliographystyle{myalpha}
\bibliography{literature}

\appendix
\numberwithin{equation}{section}
\renewcommand{\theequation}{\thesection.\arabic{equation}}
\renewcommand{\theproposition}{\thesection.\arabic{proposition}}

\section{Appendix}
\label{app:apendix}

\subsection{Proof of \texorpdfstring{\Cref{thr:visibility-circulant}}{Theorem 4.3}}\label{app:vis}
Consider the equation of motion given by \eqref{eq:f_lin_circ}.
Fix $(j,i)\in E$ such that $\|z_i(0)-z_j(0)\| = \max_{(l,k)\in E} \|z_k(0)-z_l(0)\|$. 
That is, robots $i$ and $j$ are furthest apart initially among those that do communicate. We drop the time-dependence in notation and compute
and estimate
\begin{align*}
	\frac{\mathrm{d}}{\mathrm{d}t} \|z_i-z_j\|^2 
	& = 2 \left( - \| z_i-z_j \|^2 + \sum_{k=0}^{N-1} w_k \langle z_i-z_j, z_{i+k}-z_{j+k} \rangle \right) \\
    &\le 2 \left( - \| z_i-z_j \|^2 + \sum_{k=0}^{N-1} w_k \| z_i-z_j \| \cdot \| z_{i+k}-z_{j+k} \| \right) \\
    &\le 2 \|z_i-z_j\|^2 \left(-1+\sum_{k=0}^{N-1} w_k \right) \le 0.
\end{align*}
Therein, the first estimate is due to the Cauchy-Schwarz inequality. For the second estimate, note that the sum in the third expression contains only terms $\| z_{i+k}-z_{j+k} \|$ for $(j+k,i+k)\in E$ which satisfy $\| z_{i+k}-z_{j+k} \| \le \| z_{i}-z_{j} \|$ by assumption.
In fact, since $W=(w_{i,j})_{i,j=1}^{N-1}$ is circulant, we have $w_{i+k,j+k}=w_{i,j}=w_l\ne 0$ for $j=i+l\bmod N$ and any $k \in \set{0,\dotsc,N-1}$. Thus, also $(j+k,i+k)\in E$. 
The final estimate is due to the fact that $W$ is consistent.

Hence, the distance of the two maximally distant robots cannot increase initially and by using $z_0(t),\dotsc,z_{N-1}(t)$ as a new initial configuration, the same is true for arbitrary $t>0$.

Note that we may not rule out that two robots $k,l$ with $(l,k)\in E$ that are not maximally distant increase their distance. However, assume ${\|z_k(t^*)-z_l(t^*)\|>\mathcal{C}}$ for some ${t^*>0}$. Since $\|z_i(t)-z_j(t)\|$ is non-increasing, there must be some $0<t<t^*$ at which ${\|z_k(t)-z_l(t)\|=\|z_i(t)-z_j(t)\| \le \mathcal{C}}$ by continuity of the solution $z(t)$. Then we may again use $z_0(t), \dotsc, z_{N-1}(t)$ as a new initial configuration and the argument above with the roles of $(j,i)$ and $(l,k)$ switched yields that $\|z_k(t)-z_l(t)\|$ cannot increase any further contradicting the assumption. In particular, $k$ and $l$ cannot lose sight of each other.
\qed

\subsection{Details to the Proof of \texorpdfstring{\Cref{thm:decomp}}{Theorem 4.8}}\label{app:decomp}
In general, solutions of \eqref{eq:Z_lin} are linear combinations of the so-called \emph{fundamental solutions}
\begin{align}
	\label{eq:fund_sol}
        {Z}_{j,i}^x(t)  = e^{(\lambda_j-1)t} \sum_{\ell = 0}^i \frac{t^\ell}{\ell!}(\xi_{j,\ell},0) \text{ and }
        {Z}_{j,i}^y(t)  = e^{(\lambda_j-1)t} \sum_{\ell = 0}^i \frac{t^\ell}{\ell!}(0,\xi_{j,\ell})
\end{align}
where $\lambda_0,\dotsc,\lambda_k$ are the eigenvalues of $W$ counting geometric multiplicities and 
\[
\xi_{0,1}, \dotsc, \xi_{0,m_0}, \dotsc, \xi_{k,1},\dotsc, \xi_{k,m_k} \in \C^N
\]
are the corresponding eigenvectors and generalized eigenvectors: $(W-\lambda_j\mathbf{I}_N)\xi_{j,0}=0$ and $(W-\lambda_j\mathbf{I}_N)\xi_{j,i}=\xi_{j,i-1}$ for $j=0,\dotsc,k$ and $i=2,\dotsc,m_j$. 
By \Cref{prop:circ_spectral} a circulant matrix $W$ is diagonalizable as the eigenvectors $v_j\in \C^N$ are linearly independent. 
Thus, these fundamental solutions reduce to
\begin{align*}
	{Z}_{j}^x(t) = e^{(\lambda_j-1)t}\begin{pmatrix}v_{j}\\0\end{pmatrix} \text{ and }
	{Z}_{j}^y(t) = e^{(\lambda_j-1)t}\begin{pmatrix}0\\v_{j}\end{pmatrix},
\end{align*}
i.e, $m_j = 1$ and $\xi_{j,0} = v_j$ for $j=0,\dotsc,N-1$ (cf. \eqref{eq:fund_sol}).
To construct a real-valued basis that is meaningful for interpretation in terms of configurations of $N\in \N$ robots, we proceed as follows:
First, we set 
\begin{align*}
    \T{V_0} = \Span(v_0) \subseteq \R^N, V_0 =  \Span\left(\begin{pmatrix}v_0\\0\end{pmatrix}, \begin{pmatrix}0\\ v_0\end{pmatrix}\right) \subseteq \R^{2N} \text{ and } \Lambda_0 = \lambda_0 = 1
\end{align*}
and consider two cases:

\subsubsection{\texorpdfstring{$N=2k+1$}{N=2k+1} odd:} For any $j=1,\dotsc,k$, we consider the pair $v_j \in \C^N$ and its complex conjugate $\overline{v_j}=v_{N-j}\in \C^N$ and define the subspace $\T{V}_j \subseteq \R^N$ spanned by the real- and imaginary parts of $v_j$ and $v_{N-j}$, which is
\begin{align*}
	\T{V}_j = \Span\left(\Re(v_j),\Im(v_j)\right) \subseteq \R^N.
\end{align*}
By abusing the notation and considering $\T{V}_j$ to be also the matrix containing the basis vectors of $\T{V}_j$ we immediately have

\begin{align}\label{eq:blockA}
W \T{V}_j = \T{V}_j \Lambda_j \text{ with } \Lambda_j = \begin{pmatrix}\Re(\lambda_j)& -\Im(\lambda_j)\\ \Im(\lambda_j) &\Re(\lambda_j)
\end{pmatrix}
\end{align}
since the pair of eigenvalues $\lambda_j \in \C$ and $\overline{\lambda_j} = \lambda_{N-j}$ act as a matrix $\Lambda_j \in \R^{2\times 2}$
in the subspace $\T{V}_j \subseteq \R^N$.  For brevity we will use this double meaning from now on. 
As the eigenvectors $v_j\in \C^N$ are linearly independent, we obtain a decomposition $\R^N = \bigoplus_{j=0}^{k} \T{V}_j$, in which the weight matrix $W\in \R^{N\times N}$ becomes a block-diagonal of the form
\begin{align*}
	\T{V}^{-1} W \T{V} = 
	\begin{pmatrix}
		\Lambda_0 & & & \\
		& \Lambda_1 & & \\
		& &  \ddots & \\
		& & & \Lambda_{k}
	\end{pmatrix} \in \R^{N\times N}   
\end{align*}
where $\T{V} = (\T{V}_0\dotsb \T{V}_{k}) \in \R^{N\times N}$.
Now, recall the block-diagonal structure in $\mathbf{W}\in \R^{2N\times 2N}$, which corresponds to the disconnected dynamical behavior of both coordinates of all robots. Hence, the subspace $\T{V}_j \subseteq \R^N$ can be considered as either $x$- or $y$-coordinates and we set
\begin{align*}
	V^x_j = \Span\left(\begin{pmatrix}\Re(v_j)\\0\end{pmatrix}, \begin{pmatrix}\Im(v_j)\\0\end{pmatrix}\right) \subseteq \R^{2N} \text{ and }
	V^y_j = \Span\left(\begin{pmatrix}0\\ \Re(v_j)\end{pmatrix}, \begin{pmatrix}0\\ \Im(v_j)\end{pmatrix}\right) \subseteq\R^{2N}.
\end{align*}
As both $V^x_j \subseteq \R^N$ and $V^x_j \subseteq \R^N$ correspond to the same eigenvalue, we build their direct sum and define
\begin{align}\label{eq:Vj_}
    V_j &= V_j^x \oplus V_j^y =  \Span\left(\begin{pmatrix}\Re(v_j)\\0\end{pmatrix}, \begin{pmatrix}\Im(v_j)\\0\end{pmatrix}, \begin{pmatrix}0\\ \Re(v_j)\end{pmatrix}, \begin{pmatrix}0\\ \Im(v_j)\end{pmatrix}\right).
\end{align}
which by construction is a $4$-dimensional stable invariant subspace of the model \eqref{eq:Z_lin}.
For this decomposition $\R^{2N} = \bigoplus_{j=0}^{k} V_j$ the matrix ${\mathbf{{W}}\in \R^{2N\times 2N}}$ in \eqref{eq:Z_lin} becomes
\begin{align}\label{eq:VBLOCK}
    V^{-1} \mathbf{W} V = \begin{pmatrix}
        \Lambda_0 & & & \\
        & \Lambda_0 & & \\
        & & \ddots & \\
        & & & \Lambda_{k}\\
        & & & & \Lambda_{k}
    \end{pmatrix} \in \R^{2N\times 2N}
\end{align}
where $V = \left(V_0 \dotsb V_{k}\right) \in \R^{2N\times 2N}$. By definition of $V_j\subseteq \R^{2N}$ in \eqref{eq:Vj_} the spanning vectors correspond to configurations, where all robots have only non-zero components in one coordinate ($x$ or $y$). Thus, in particular for visualization purposes of the configurations contained in $V_j\subseteq\R^{2N}$ in the Euclidean plane, we propose the following basis instead
\begin{align}\label{eq:VjA}
	V_j = \Span\left(
	\begin{pmatrix}\Re(v_j)\\ -\Im(v_j)\end{pmatrix}, \begin{pmatrix}\Im(v_j)\\\Re(v_j)\end{pmatrix}, \begin{pmatrix}-\Im(v_j)\\ \Re(v_j)\end{pmatrix},\begin{pmatrix}\Re(v_j)\\ \Im(v_j)\end{pmatrix}
	\right),
\end{align}
(cf. \eqref{eq:Vj}). Note that, by using this explicit adapted basis, the block structure into blocks $\Lambda_j$ in \eqref{eq:VBLOCK} is not changed. In fact, the change-of-base matrix is given by
\begin{align*}
	T = \frac{1}{2}\begin{pmatrix}
		1&0&0&-1\\
		0&1&1&0\\
		0&-1&1&0\\
		1&0&0&1
	\end{pmatrix}
	\quad \text{resp.}\quad
	T^{-1} = \begin{pmatrix}
		1&0&0&1\\
		0&1&-1&0\\
		0&1&1&0\\
		-1&0&0&1
	\end{pmatrix}
\end{align*}
which, after changing to the adapted basis \eqref{eq:VjA}, yields the same block matrices as
\begin{align*}
	T \begin{pmatrix}
		\Lambda_j&\\
		&\Lambda_j
	\end{pmatrix}
	T^{-1} = \begin{pmatrix}
		\Lambda_j&\\
		&\Lambda_j
	\end{pmatrix}
\end{align*}
In particular, the doubled blocks are not merged into $4 \times4$-dimensional blocks.

\subsubsection{\texorpdfstring{$N = 2k$}{2k} even:} Here, we proceed for $j=1,\dotsc,k-1$ the same way as in the odd case. Since the eigenvalue $\lambda_k \in \R$ is real (cf. \Cref{rem:circ_EW} (b)) we can simply define
\begin{align*}
    \T{V}_k = \Span(v_k) \subseteq \R^N \text{ and } \Lambda_k = \lambda_k \in \R.
\end{align*}
Here, \eqref{eq:blockA} reduces to the eigenvalue equation $Wv_k = \Lambda_k v_k$ for $\Lambda_k = \lambda_k \in \R$. The spaces $V_k^x$ and $V_j^y$ are given by
\begin{align*}
    V^x_k = \Span\left(\begin{pmatrix}v_k\\0\end{pmatrix}\right) \subseteq \R^{2N} \text{ and }
    V^y_k = \Span\left(\begin{pmatrix}0\\ v_k\end{pmatrix}\right) \subseteq\R^{2N}.
\end{align*}
and their direct sum is
\begin{align*}
    V_k = V_k^x \oplus V_k^y = \Span\left(\begin{pmatrix}v_k\\0\end{pmatrix},\begin{pmatrix}0\\ v_k\end{pmatrix}\right),
\end{align*}
which by construction is a $2$-dimensional stable invariant subspace of \eqref{eq:Z_lin}. Note that the block-structure of \eqref{eq:VBLOCK} still holds true for this case.

\subsection{Details to the Proof of \texorpdfstring{\Cref{thm:decomp_nonlinear}}{Theorem 5.1}}\label{app:decomp_nonlinear}
\subsubsection{Representation of \texorpdfstring{$\normi$}{normi}:}
It turns out, that $\normi$ has a useful representation on each of the subspaces $V_j^\pm$. For $a,b\in \R$ let
\begin{align*}
	{Z} = V_j^i \begin{pmatrix}
		a\\b
	\end{pmatrix} = \begin{pmatrix}\Re(v_j)& \pm\Im(v_j)\\ \mp\Im(v_j) & \Re(v_j)\end{pmatrix} \begin{pmatrix}
		a\\b
	\end{pmatrix}
	= a\begin{pmatrix}\Re(v_j)\\ \mp\Im(v_j)\end{pmatrix} + b \begin{pmatrix}\pm\Im(v_j)\\\Re(v_j)\end{pmatrix}
\end{align*}
and we compute
\begin{align*}
	\normi(Z) &= \normi\left(a \begin{pmatrix}\Re(v_j)\\ \mp\Im(v_j)\end{pmatrix} + b \begin{pmatrix}\pm\Im(v_j)\\\Re(v_j)\end{pmatrix}\right)\\
		& = \norm{(a\Re(v_j)_i\pm b\Im(v_j)_i,\mp a\Im(v_j)_i+b\Re(v_j)_i)}\\
	    & = \sqrt{(a\Re(v_j)_i\pm b\Im(v_j)_i)^2 + (\mp a\Im(v_j)_i+b\Re(v_j)_i)^2}\\
	    & = \sqrt{(a^2+b^2)(\Re(v_j)_i^2 + \Im(v_j)_i^2)}\\
	    & = \sqrt{a^2+b^2}\\
	    & = \norm{(a,b)}
\end{align*}
which only depends the coefficients $(a,b)\in \R^2$ and, in particular, each entry $i$ is the same. The step from the fourth to the fifth line holds true, since all entries of $v_j$ are primitive roots of unity.

\subsubsection{Derivation of \texorpdfstring{\eqref{eq:dotnorm}}{(19)}:}
Using \eqref{eq:ab} and the chain rule we compute for $(a,b)\in \R^2$:
\begin{align*}
	\frac{\mathrm{d}}{\mathrm{d} t}\norm{\begin{pmatrix}a\\b\end{pmatrix}} (t) & = \mathrm{D}\norm{\cdot}(a,b) \cdot \frac{\mathrm{d}}{\mathrm{d}t}{\begin{pmatrix}a\\b\end{pmatrix}}\\
		& = \frac{1}{\norm{(a,b)}} (a,b) \n(\abs{-1 + \lambda_j}\norm{(a,b)}) (-\mathbf{I}_2 + \Lambda_j )\begin{pmatrix}a\\b\end{pmatrix}\\
		& = \n(\abs{-1 + \lambda_j}\norm{(a,b)}) \frac{(a,b)
			\begin{pmatrix}	-1+\Re(\lambda_j) & \Im(\lambda_j)\\-\Im(\lambda_j) & -1+\Re(\lambda_j)\end{pmatrix}
			\begin{pmatrix}a\\b\end{pmatrix}}
			{\norm{(a,b)}}\\
		& = \n(\abs{-1 + \lambda_j}\norm{(a,b)}) \frac{(-1+\Re(\lambda_j))(a^2 + b^2)}{\norm{(a,b)}}\\
		& =  (-1+\Re(\lambda_j)) \n(\abs{-1 + \lambda_j}\norm{(a,b)}) \norm{\begin{pmatrix}a\\b\end{pmatrix}} < 0.
\end{align*}

\subsection{Proof of \texorpdfstring{\Cref{prop:visibility-circulant-nonlinear}}{Proposition 5.2}}\label{app:vis_nonlinear}
For $Z \in V_j^\pm$ recall the equation of motion \eqref{eq:DS_TZ} which can be written as
\begin{align*}
	\dot{z}_i & = \n(-1 + \abs{\lambda_j}\norm{(a,b)}) \left(-z_i + \sum_{\ell=0}^{N-1} w_\ell z_{i+\ell}\right)
\end{align*}
Here, the scaling factor $\n(-1 + \abs{\lambda_j}\norm{(a,b)})$ is the same for each robot $i$. This allows us to proceed as in the linear case (cf. \Cref{app:vis}).

Fix $(q,p)\in E$ such that $\|z_p(0)-z_q(0)\| = \max_{(l,k)\in E} \|z_k(0)-z_l(0)\|$ and we estimate as before
\begin{align*}
	\frac{\mathrm{d}}{\mathrm{d}t} \|z_p-z_q\|^2 
		& = 2\n(-1 + \abs{\lambda_j}\norm{(a,b)})\left( - \| z_p-z_q \|^2 + \sum_{k=0}^{N-1} w_k \langle z_p-z_q, z_{p+k}-z_{q+k} \rangle \right) \\
		&\le 2 \n(-1 + \abs{\lambda_j}\norm{(a,b)})\left( - \| z_p-z_q \|^2 + \sum_{k=0}^{N-1} w_k \| z_p-z_q \| \cdot \| z_{p+k}-z_{q+k} \| \right) \\
		&\le 2 \n(-1 + \abs{\lambda_j}\norm{(a,b)})\|z_p-z_q\|^2 \left(-1+\sum_{k=0}^{N-1} w_k \right) \le 0.
\end{align*}
Therein, we exploited the fact that the nonlinear scaling is independent of the robots' index. For the rest of the proof we refer to argumentation in \Cref{app:vis}.
\qed

\subsection{Proof of \texorpdfstring{\Cref{prop:finite_time}}{Proposition 5.3}}\label{app:finite_time}
For $\n(\norm{(x,y)}) = \norm{(x,y)}^{-1}$, \eqref{eq:dotnorm} can be simplified to
\begin{align*}
	\dot{\norm{\begin{pmatrix}a\\b\end{pmatrix}}} = (-1+\Re(\lambda_j))\frac{\norm{(a,b)}}{\abs{-1+\lambda_{j}}\norm{(a,b)}} = \frac{-1+\Re(\lambda_j)}{\abs{-1+\lambda_{j}}} = \mbox{const}\leq 0.
\end{align*} 
The solution of this ordinary differential equation is given by
\begin{align*}
	\norm{\begin{pmatrix}a\\b\end{pmatrix}}(t) = \norm{\begin{pmatrix}	a_0\\b_0\end{pmatrix}} + \frac{-1+\Re(\lambda_{j})}{\abs{-1+\lambda_j}}t,
\end{align*}
i.e, the decay of $\norm{(a,b)}$ is \emph{linear} with rate $\frac{-1+\Re(\lambda_j)}{\abs{-1+\lambda_j}}$ until $\norm{(a,b)} = 0$. In particular, gathering is reached in finite time
\begin{align*}
	t_{\mbox{gather},j} = -\frac{\abs{-1+\lambda_{j}}}{-1+\Re(\lambda_j)}\norm{\begin{pmatrix}a_0\\b_0\end{pmatrix}} = -\frac{\abs{-1+\lambda_{j}}}{-1+\Re(\lambda_j)} \normi(Z_0) > 0\quad \text{ for any } i=1,\dotsc,N.
\end{align*}
\qed

\subsection{Examples}\label{app:example}
\subsubsection{\gtm{} Protocol}
As a second example we consider the \gtm{} protocol (cf. \cite{Kling.2011}). In contrast to the \nbug{} problem, the underlying interaction graph is symmetric. Its weight matrix is given by ${W = \circulant(0,\frac{1}{2},0,\dotsc,\frac{1}{2})}$, i.e., robot $i$ is influenced by its first neighbors to the left and the right. Its symmetric interaction graph is illustrated in \Cref{fig:sym_graph-a} and \Cref{thr:gathering-circulant} implies that the only gathering linear protocol on such a undirected next neighbor graph is the \gtm{} protocol.

\begin{figure}[!htb]
    \centering
    \hfil
    \subfloat[][\centering \gtm{} protocol]{
        \label{fig:sym_graph-a}
        \resizebox{!}{0.25\textwidth}{
            \centering
\begin{tikzpicture}[<->,
	>=stealth',
	shorten >=1pt,
	auto,
	node distance=1cm,
	main node/.style={line width=1.5pt, circle, scale = 2, draw, font=\sffamily\tiny, inner sep=1pt}]
	\def\ngon{7}
	\node[main node, scale = 1.5, regular polygon, regular polygon sides=\ngon, minimum size=1cm, draw=none] (p) {};
	\node[font=\sffamily\footnotesize, inner sep=1pt, scale=1.5](d) at (p.corner 1) {...};
	\node[main node, minimum size=.3cm](N) at (p.corner 2) {N};
	\node[main node, minimum size=.3cm](1) at (p.corner 3) {1};
	\node[main node, minimum size=.3cm](2) at (p.corner 4) {2};
	\node[main node, minimum size=.3cm](3) at (p.corner 5) {3};
	\node[main node, minimum size=.3cm](4) at (p.corner 6) {4};
	\node[main node, minimum size=.3cm](5) at (p.corner 7) {5};
	
%	\draw[red, opacity=.3] (0,0) circle (2.25cm);
	
	\path[every node/.style={font=\sffamily\small}, line width =1pt]
    
	(1) edge [bend left=10] node {} (N)
	(2) edge [bend left=10] node {} (1)
	(3) edge [bend left=10] node {} (2)
	(4) edge [bend left=10] node {} (3)
	(5) edge [bend left=10] node {} (4)
	(d) edge [bend left=10] node {} (5)
	(N) edge [bend left=10] node {} (d)
	
	;
\end{tikzpicture}%
        }
    }
    \hfil
    \subfloat[][\centering \gta{} protocol]{
        \label{fig:sym_graph-b}
        \resizebox{!}{0.25\textwidth}{
            \centering
\begin{tikzpicture}[<->,
	>=stealth',
	shorten >=1pt,
	auto,
	node distance=1cm,
	main node/.style={line width=1.5pt, circle, scale = 2, draw, font=\sffamily\tiny, inner sep=1pt}]
	\def\ngon{7}
	\node[main node, scale = 1.5, regular polygon, regular polygon sides=\ngon, minimum size=1.2cm, draw=none] (p) {};
	\node[font=\sffamily\footnotesize, inner sep=1pt, scale=1.5](d) at (p.corner 1) {...};
	\node[main node, minimum size=.3cm](N) at (p.corner 2) {N};
	\node[main node, minimum size=.3cm](1) at (p.corner 3) {1};
	\node[main node, minimum size=.3cm](2) at (p.corner 4) {2};
	\node[main node, minimum size=.3cm](3) at (p.corner 5) {3};
	\node[main node, minimum size=.3cm](4) at (p.corner 6) {4};
	\node[main node, minimum size=.3cm](5) at (p.corner 7) {5};
	
%	\draw[red, opacity=.3] (0,0) circle (2.25cm);
	
	\path[every node/.style={font=\sffamily\small}, line width =1pt]
	% (d) edge [] node {} (N)
	(d) edge [] node {} (1)
 	(d) edge [] node {} (2)
	(d) edge [] node {} (3)
	(d) edge [] node {} (4)
	% (d) edge [] node {} (5)

    (N) edge [] node {} (2)
	(N) edge [] node {} (3)
 	(N) edge [] node {} (4)
	(N) edge [] node {} (5)

 	(1) edge [] node {} (3)
	(1) edge [] node {} (4)
	(1) edge [] node {} (5)

	(2) edge [] node {} (4)
	(2) edge [] node {} (5)

	(3) edge [] node {} (5)
    
	(1) edge [bend left=10] node {} (N)
	(2) edge [bend left=10] node {} (1)
	(3) edge [bend left=10] node {} (2)
	(4) edge [bend left=10] node {} (3)
	(5) edge [bend left=10] node {} (4)
	(d) edge [bend left=10] node {} (5)
	(N) edge [bend left=10] node {} (d)

    (1) edge [loop, out=200,in=150,looseness=6] node {} (1)
	(2) edge [loop, out=250,in=200,looseness=6] node {} (2)
 	(3) edge [loop, out=330,in=280,looseness=6] node {} (3)
 	(4) edge [loop, out=380,in=330,looseness=6] node {} (4)
   	(5) edge [loop, out=430,in=380,looseness=6] node {} (5)
    (N) edge [loop, out=150,in=100,looseness=6] node {} (N)
    (d) edge [loop, out=480,in=430,looseness=25] node {} (d)

	;
\end{tikzpicture}%
        }
    }
    \caption{Illustration of the circulant interaction graphs of the \gtm{} and \gta{} protocols. The arrangement of the vertices does not reflect the physical location of robots.}
\end{figure}

Since the \gtm{} protocol is symmetric, we have real-valued eigenvalues $\lambda_j \in \R$ (cf. \Cref{rem:circ_EW} (b)), i.e., they coincide with the convergence rates $\Re(\lambda_j)$. Using the eigenvalue formula \eqref{eq:EW} with ${w=(0,\frac{1}{2},0,\dotsc,\frac{1}{2})^T \in \R^N}$ we compute
\[ 
\lambda_j = \Re(\lambda_j) = \cos\left(\frac{2\pi j}{N}\right) \text{ for } j=0,\dots, N-1.
\]
Observe that the convergence rates for this protocol are the same as for the \nbug{} problem discussed in \Cref{ex:EW} whose rates are illustrated in \Cref{fig:conv_rates}.
In particular, we obtain the same hierarchical decomposition. Note that, from the dynamical systems perspective both models gather at the same speed. As the decomposition \eqref{eq:decomp} is independent of the generating vector $w\in \R^N$, the generating configurations are also visualized in \Cref{fig:EV_C}.  

\subsubsection{\gta{} Protocol}
Finally as a third example, we take a look at the (global) \gta{} protocol (also called \gtg{} in \cite{Cohen.2004b}). In this model every robot has global vision and is therefore influenced by all robots (including itself). Hence, the corresponding interaction graph is complete and its weight matrix can be written as $W = \circulant(\frac{1}{N},\dotsc,\frac{1}{N})$. In \Cref{fig:sym_graph-b} its complete interaction graph is illustrated. 
Similarly, \Cref{thr:gathering-circulant} implies that the \gta{} protocol is gathering. It is, in fact, the only gathering linear protocol on the complete undirected graph for which all weights are equal.

The generating vector $w=(\frac{1}{N},\dotsc,\frac{1}{N})^T \in \R^N$ of the \gta{} protocol yields the eigenvalues
\begin{align*}
    \lambda_j = \frac{1}{N} \sum_{i=0}^{N-1}\omega^{ij} = \begin{cases}
        1, &\text{if } j=0,\\
        0, &\text{if } j\neq 0,
    \end{cases}
\end{align*}
since the average of the $N$-th roots of unity vanishes for $j\neq0$. Again, by the symmetry of the protocol all eigenvalues $\lambda_j$ are indeed real and thus coincide with the convergence rates. In particular, for this protocol we have $\Re(\lambda_j) = 0$ for all $j\neq 0$ and every configuration converges with the same speed. Even though \Cref{thm:decomp} decomposes the state space into stable subspaces $V_j\subseteq \R^{2N}$ (also shown in \Cref{fig:EV_C}), they cannot be distinguished by their convergence rates $\Re(\lambda_j)\in \R$.

\end{document}